\documentclass[a4paper,11pt]{amsart}
\usepackage[top=2.2cm,left=2.8cm,right=2.8cm,bottom=2.8cm]{geometry}
\usepackage[foot]{amsaddr}

\usepackage{hyperref}
\hypersetup{
	colorlinks,
	citecolor=green,
	filecolor=black,
	linkcolor=blue,
	urlcolor=blue
}
\hypersetup{linktocpage}

\usepackage{cite}
\usepackage{graphicx}
\usepackage{amsbsy}
\usepackage{latexsym}
\usepackage{amsfonts}
\usepackage{amssymb}
\usepackage{upgreek}
\usepackage{mathtools}
\usepackage{cleveref}
\usepackage[usenames]{color}
\usepackage{amsmath,amsthm}
\usepackage{enumerate}
\usepackage{mathrsfs}
\usepackage{stmaryrd}
\usepackage{dsfont}
\usepackage{varwidth}
\usepackage{pdfpages}
\usepackage{comment}

\usepackage{tikz}
\usetikzlibrary{arrows.meta,calc,positioning,fit}
\usetikzlibrary{decorations.pathreplacing}

\newcommand{\bG}{\ensuremath{\mathbb{G}}}

\newcommand{\bN}{\ensuremath{\mathbb{N}}}

\newcommand{\bR}{\ensuremath{\mathbb{R}}}

\newcommand{\bW}{\ensuremath{\mathbb{W}}}

\newcommand{\cR}{\ensuremath{\mathcal{R}}}

\newcommand{\N}{\bN}
\newcommand{\R}{\bR}

\renewcommand{\d}{\mathrm{d}}

\newcommand{\dx}{\mathrm{d}x}

\newcommand{\id}{I}

\DeclareMathOperator{\Span}{span}

\DeclareMathOperator*{\argmin}{arg\,min}
\DeclareMathOperator*{\argmax}{arg\,max}

\newcommand{\df}{\mathcal A}
\newcommand{\adv}{\mathcal B}
\newcommand{\fp}{\mathrm{fp}}

\newcommand{\degone}{\mathrm{1}}
\newcommand{\deghigher}{>}

\newcommand{\D}{{\mathrm D}}

\providecommand{\abs}[1]{\lvert#1\rvert}
\providecommand{\bigabs}[1]{\bigl\lvert#1\bigr\rvert}

\providecommand{\biggabs}[1]{\biggl\lvert#1\biggr\rvert}

\providecommand{\norm}[1]{\lVert#1\rVert}
\DeclarePairedDelimiter{\inner}{\langle}{\rangle}
\providecommand{\bignorm}[1]{\bigl\lVert#1\bigr\rVert}
\providecommand{\Bignorm}[1]{\Bigl\lVert#1\Bigr\rVert}

\DeclarePairedDelimiter{\set}{\{}{\}}

\newtheorem{theorem}{Theorem}
\newtheorem{lemma}[theorem]{Lemma}
\newtheorem{prop}[theorem]{Proposition}
\newtheorem{cor}[theorem]{Corollary}

\theoremstyle{definition}

\newtheorem{example}[theorem]{Example}

\theoremstyle{remark}
\newtheorem{remark}[theorem]{Remark}

\numberwithin{equation}{section}
\numberwithin{theorem}{section}

\theoremstyle{plain}
\newtheorem{assumption}[theorem]{Assumption}

\newcommand{\Chi}{\raise .3ex
	\hbox{\large $\chi$}}

\newcommand{\cond}{\, : \,}

\usepackage[section]{algorithm}
\usepackage{algorithmicx}
\usepackage{algpseudocode}

\usepackage{todonotes}

\title[Adaptive Gaussian Splitting Scheme]{An Adaptive Gaussian Splitting Scheme for Linear Fokker--Planck Equations}
\makeatletter\def\@setemails{\mbox{{\itshape E-mail addresses}:\space}{\ttfamily\emails}.}\makeatother
\author{Markus Bachmayr$^{1}$, Tianyu Jin$^{1}$, Damiano Lombardi$^{2}$, Olga Mula$^{3}$ and Ahmad Zorkot$^{3}$}
\address{$^{1}$Institut f\"ur Geometrie und Praktische Mathematik, RWTH Aachen University, Im S\"usterfeld 2, 52072 Aachen, Germany}
\address{$^{2}$Laboratoire Jacques-Louis Lions, Sorbonne Universit\'e, and Inria Paris, 48 rue Barrault, 75013, Paris, France}
\address{$^{3}$Faculty of Mathematics, University of Vienna, Oskar-Morgenstern-Platz 1, 1090 Vienna, Austria}
\email{bachmayr@igpm.rwth-aachen.de}
\email{jin@igpm.rwth-aachen.de}
\email{damiano.lombardi@inria.fr}
\email{olga.mula.hernandez@univie.ac.at}
\email{ahmad.zorkot@univie.ac.at}

\thanks{M.B. and T.J.\ acknowledge funding by the European Union (ERC, COCOA, 101170147).}

\begin{document}

\maketitle

\begin{abstract}
  We develop an adaptive numerical method for linear Fokker--Planck equations on $\bR^d$, based on approximating the evolving probability density by linear combinations of Gaussian functions (LCGs). The method applies to polynomial-type drift fields, without requiring a gradient-flow structure or a stationary distribution. It relies on an operator splitting that separates diffusion and affine drift, for which Gaussian functions can be propagated exactly, from higher-order transport, whose evolution is approximated adaptively by LCGs. Both the number of Gaussian components and their means and covariances are adjusted to meet the prescribed accuracy tolerances. The adaptivity is driven by computable error indicators involving polynomial moments of Gaussians that can be evaluated analytically. We provide a rigorous error analysis that relates the approximation error to the indicators driving adaptivity and illustrate the performance of the scheme in several test cases.

  \smallskip

\end{abstract}

\section{Introduction}\label{sec:introduction}
Fokker--Planck equations describe the evolution of probability densities under deterministic drift and stochastic fluctuations, with applications in molecular simulations, uncertainty quantification, and machine learning. When the joint probability distribution of a large number of correlated random variables is described in this manner, this leads to a time-dependent partial differential equation on a high-dimensional domain.
We consider problems of this type on $\bR^d$, for which mesh-free methods are particularly attractive. The method we consider is based on approximations by linear combinations of Gaussian functions (LCGs), where both the number and the parameters of the Gaussian terms are adjusted during the evolution. In particular, the adaptation of covariance parameters of the individual Gaussians enables the efficient approximation of anisotropic structures without preferred directions.

More specifically, we consider the linear Fokker--Planck equation
\begin{equation}\label{eq:FP}
  \begin{aligned}
    \partial_t u(t,x) & = \kappa\Delta u(t,x) - \nabla\cdot\bigl(F(x)u(t,x)\bigr), & \qquad & (t,x)\in(0,T]\times\R^d, \\
    u(0,x)            & = u_0(x),                                                  &        & x\in\R^d,
  \end{aligned}
\end{equation}
where $T>0$, $\kappa>0$ is the diffusion coefficient, and $F:\R^d\to\R^d$ is the drift field. We focus on polynomial-type fields, decomposed as $F=F_{\deghigher}+F_{\degone}$ into an affine part $F_{\degone}$ and higher-order terms $F_{\deghigher}$. Our method splits the evolution into diffusion with affine drift $F_{\degone}$ and transport driven by $F_{\deghigher}$. The first subflow propagates every Gaussian exactly, with explicit formulas for its mean and covariance (see Lemma~\ref{lem:affine_gaussian_preservation}). The second generally does not leave the class of LCGs invariant, so its output is approximated again by an LCG. The number of components, their means, and their covariances are selected adaptively using computable error indicators. Crucially, the quantities entering these indicators reduce to polynomial moments of Gaussians and can be evaluated analytically. The resulting adaptation operates directly on $\bR^d$ and requires neither a spatial grid nor sampling. It does not require $F$ to be conservative or to admit a stationary distribution either. Efficient selection of the moments needed in high dimensions is an important further question, which we leave for future work. Here, we confine the discussion to doing adaptivity with LCGs and reaching a target tolerance.

There exists an extensive literature on numerical schemes for problems of the basic type \eqref{eq:FP}. Classical discretizations include finite difference schemes \cite{CHANG19701}, finite element methods \cite{masud_bergman2005}, and spectral or spectral element methods \cite{lozinski2003fast,KS2009}. Some methods are designed to preserve positivity or mass conservation \cite{CHANG19701}. Adaptive spatial discretizations have been developed for time-dependent Fokker--Planck problems \cite{MW2018} and related kinetic models \cite{RHCTT2024}. Such methods typically rely on a spatial discretization and are posed on bounded domains. When applied to whole-space problems, they may require domain truncation. Particle methods avoid a spatial mesh but typically use stochastic trajectories and statistical density reconstruction \cite{DaumHuang2003,Robert2004}, where reconstruction of the full probability density may require a prohibitively large number of samples even in moderate dimensions. Deterministic particle variants have also been developed for Fokker--Planck equations with a suitable gradient flow structure. For instance, the blob method \cite{CarrilloCraigPatacchini2019} represents diffusion through regularized particle interactions. In place of a mesh, such methods still require a sufficiently rich particle representation of the underlying probability distribution.

A different mesh-free strategy represents the solution through time-dependent parameters, as in neural parametric and Neural Galerkin methods \cite{liu_li_zha_zhou2022, BCM2025}, or other types of parametric functions, including linear combinations of Gaussians \cite{LLB2026}. These methods commonly use sampling to approximate the quantities governing parameter evolution, and dynamical sampling methods need to be employed \cite{wen2023couplingparameterparticledynamics,Bruna2024,BCM2025}. Gaussian parametrizations have a long history in time-dependent wave and quantum problems \cite{LasserLubich2020,JinMarkowichSparber2011}, and several Gaussian-based solvers have been proposed for Fokker--Planck dynamics. Fixed-size Gaussian representations are considered in \cite{TerejanuSinglaSinghScott2008,AndersonFarazmand2024, LLB2026}, while adaptive Gaussian approximations, including residual-based control and component splitting or merging, appear in \cite{scharpenberg_lukacova2012,MukherjeeRaiSinglaSinghPatra2016,vishwajeet_singla2018}. Other approaches exploit special conditional-Gaussian structure \cite{chen_majda2018} or sum-of-Gaussians representations of kernels for fractional Fokker--Planck equations \cite{jiang_wang_zhou2026}. Gaussian approximation and other types of parametric approximation have also been considered when viewing the Fokker--Planck equation as a Wasserstein gradient flow (see, e.g., \cite{caucheteux2026unifying, gonnella2026dynamical}), or in other metric spaces such as Bures--Wasserstein and Hellinger--Kantorovich (see \cite{LambertChewiBachBonnabelRigollet2022, liero2025evolution}). However, adaptivity aspects are not present in these works.

The distinguishing feature of our approach is a splitting between exact Gaussian propagation for the diffusive subflow with affine drift and a deterministic, moment-based adaptivity for the transport subflow of higher-order drift terms. The latter step can increase the number of Gaussian components based on computable indicators; this potential growth is controlled by an additional recompression of LCGs up to controlled errors in each time step.
The scheme does not require conservative forces or a steady state. We provide a rigorous error analysis, in particular establishing controlled errors for each step of the adaptive scheme.

The remainder of the paper is organized as follows. Section~\ref{sec:wellposedness} describes the LCG class of functions, and the functional framework in which the exact evolution takes place. Section~\ref{sec:main_scheme} introduces the splitting framework, the exact Gaussian update for the affine Fokker--Planck subproblem, and the approximated integration of the transport subproblem with nonlinear drift. Section~\ref{sec:adaptive} develops the mass-constrained optimization used for Gaussian fitting in the transport step, including its initialization and gradient evaluation. It also describes how we dynamically increase or decrease the number of Gaussians in our approximation. Section~\ref{sec:error_analysis} establishes error estimates for the proposed methods in weighted $L^2$ spaces. Numerical results for four test problems are presented in Section~\ref{sec:numerical_examples}, and concluding remarks are given in Section~\ref{sec:conclusion}.

\section{Gaussian approximations and exact evolution}\label{sec:wellposedness}

We aim to approximate solutions of Fokker--Planck equations by finite linear combinations of Gaussian functions (LCGs). We first describe in detail this class of functions from which the approximations are selected, then consider the properties of the semigroup associated to the Fokker--Planck equation under our assumptions.

For $x\in\R^d$, $|x|$ denotes the Euclidean norm, and we abbreviate $L^p(\R^d)$ by $L^p$. For a positive weight $\rho$, we denote by $L^2(\rho)$ the Hilbert space with measure $\rho(x)\,\d x$, so that $L^2(1)=L^2$.
For $s \in \N$, we write $H^s(\rho)$ for the weighted Sobolev space with norm
\[
  \norm{v}_{H^s(\rho)}^2
  =\sum_{|\alpha|\leq s} \norm{ \partial^\alpha v }_{L^2(\rho)}^2 \vspace{-6pt}
\]
and abbreviate $H^s=H^s(1)$.

\subsection{Approximation with linear combinations of Gaussians}

We represent Gaussian functions using a Cholesky-type factorization for the representation of symmetric positive definite matrices: let
\[
  \mathbb{U} = \bigl\{ U\in\bR^{d\times d} \colon U \text{ is upper triangular and } \det(U)\neq 0 \bigr\}.
\]
For $m\in\bR^d$ and $U\in\mathbb{U}$, we define the associated Gaussian function as
\begin{equation*}
  g_{m,U}(x) = \exp\left( -\frac12(x-m)^\top U^\top U(x-m) \right), \qquad x\in\bR^d.
\end{equation*}
The maximizer $m$ is also called its \emph{center;} it coincides with the expectation of the Gaussian random variable with distribution given by normalization of $g_{m,U}$.
For this distribution, $P = U^\top U$ is called the \emph{precision matrix} (whose inverse is the corresponding covariance matrix), and we adopt this terminology also for Gaussian functions without normalization. We denote the dictionary of Gaussian functions with arbitrary positive definite precision matrices by
\[
  \mathbb{G} = \bigl\{ g_{m,U} \colon m\in\bR^d,\; U\in\mathbb{U} \bigr\}.
\]
For a linear combination of $k$ Gaussian components, it will be convenient to separate the linear coefficients from the nonlinear Gaussian parameters. To this end, let $\mathbb{W} = \bR^d\times\mathbb{U}$.
For $w\in\mathbb{W}^k$, i.e. $w_i=(m_i,U_i)\in\mathbb{W}, i=1,\dots,k$, the associated Gaussian function is denoted by
\[
  g_{w_i}(x)=g_{m_i,U_i}(x).
\]
Given a coefficient vector $c=(c_1,\dots,c_k)^\top\in\bR^k$, the corresponding $k$-term linear combination of Gaussians is then written as
\begin{equation}\label{eq:lcg}
  u_{c,w}(x) = \sum_{i=1}^k c_i g_{w_i}(x),
\end{equation}
and we denote by
\begin{equation*}
  \bG_k = \set{ u_{c, w} \cond (c, w)\in \bR^k\times \bW^k }
\end{equation*}
the set of all possible linear combinations of $k$ Gaussian functions.

In our numerical scheme, the solution is approximated at time $t$ as
\[
  u(t,\cdot) \approx u_{c(t), w(t)} = \sum_{i=1}^{k(t)} c_i(t)g_{w_i(t)}\; \in \bG_{k(t)}, \qquad c(t)\in\R^{k(t)}\text{ and }w(t)\in\mathbb W^{k(t)},
\]
where both the number $k(t)$ of Gaussian components and their parameters will be chosen adaptively in the evolution.

We conclude this section with a norm equivalence between $L^2$
and $L^2(\rho_a)$ on finite-dimensional spaces spanned by Gaussians, where
\begin{equation}
  \label{eq:rhoa}
  \rho_a(x) = \exp(a \abs{x}^2),\quad a\geq 0.
\end{equation}
The following lemma identifies the admissible values of $a$
and guides our choice of spaces for the exact solution in the
next section.

\begin{lemma}
  \label[lemma]{lemma:norm-equivalence-guassian-subspaces}
  Let $V_k = \operatorname{span}\set{g_{w_i}}_{i=1}^k$, where the
  Gaussian functions $g_{w_i}$ with $w_i = (m_i, U_i)$ are linearly independent, and define
  \begin{equation}
    \label{eq:astar}
    a^* =
    \min_{1\leq i\leq k}\lambda_{\min}(U_i^\top U_i).
  \end{equation}
  For $a\geq 0$, the norms $\norm{\cdot}_{L^2}$ and
  $\norm{\cdot}_{L^2(\rho_a)}$ are equivalent on $V_k$
  if and only if $0\leq a<a^*$, in which case for all $u\in V_k$,
  \begin{equation}\label{eq:unormeq}
    \norm{u}_{L^2}
    \leq \norm{u}_{L^2(\rho_a)}
    \leq C(a,V_k)\norm{u}_{L^2}.
  \end{equation}
  The smallest possible constant in the upper bound is given by
  \begin{equation}
    \label{eq:upper-constant-norm-eq}
    C(a,V_k)
    =
    \sqrt{\lambda_{\max}
      \left(G_0^{-1/2}G_aG_0^{-1/2}\right)},
  \end{equation}
  with the Gramian $G_s = \bigl(\inner{g_{w_i},g_{w_j}}_{L^2(\rho_s)}\bigr)_{i,j=1}^k \in\bR^{k\times k}$ for $0\leq s<a^*$.
\end{lemma}

\begin{proof}
  To show the equivalence, we first assume $0\leq a<a^*$. For each $j$, the matrix
  $Q_j=P_j-aI_d$ with $P_j=U_j^\top U_j$ is positive definite. With $0\leq a < a^*$, completing the square
  gives
  \begin{align*}
    \norm{g_{w_j}}_{L^2(\rho_a)}^2
     & =
    \int_{\bR^d}
    \exp\left(
    -x^\top Q_jx+2x^\top P_jm_j-m_j^\top P_jm_j
    \right)\,dx \\
     & =
    \frac{\pi^{d/2}}{\sqrt{\det Q_j}}
    \exp\left(
    m_j^\top P_jQ_j^{-1}P_jm_j-m_j^\top P_jm_j
    \right)
    <\infty .
  \end{align*}
  Hence $V_k \subset L^2(\rho_a)$. Since $g_{w_j} \in L^2$, we also have that $V_k \subset L^2$. We next show that both norms are equivalent on $V_k$. For any $u=\sum_{i=1}^k c_i g_{w_i} \in V_k$, we have $\norm{u}_{L^2}^2=c^\top G_0c$ and $\norm{u}_{L^2(\rho_a)}^2=c^\top G_a c$.
  With $z=G_0^{1/2}c$, it follows that
  \begin{equation*}
    \norm{u}_{L^2(\rho_a)}^2
    =z^\top G_0^{-1/2}G_aG_0^{-1/2}z
    \leq
    \lambda_{\max}\left(G_0^{-1/2}G_aG_0^{-1/2}\right)
    \norm{u}_{L^2}^2.
  \end{equation*}
  Since $\rho_a(x)\geq1$, we conclude \eqref{eq:unormeq} with $C(a,V_k)$ as in \eqref{eq:upper-constant-norm-eq}.
  The constant is optimal by the Rayleigh quotient
  characterization of the largest eigenvalue.

  To conclude the equivalence, suppose finally that $a\geq a^*$. It then suffices to observe that we can choose $j$ such that $\lambda_{\min}(U_j^\top U_j) = a^*$, which entails $g_{w_j} \notin L^2(\rho_a)$ while $g_{w_j} \in L^2$.
\end{proof}

\subsection{Boundedness of flows}

In what follows, we consider the boundedness properties of the flow $\phi_t^{\fp}$ associated to the Fokker--Planck equation \eqref{eq:FP} for $t\geq 0$. We focus on the case of polynomial $F$ and initial data $u_0$ that can be written as arbitrary finite linear combinations of Gaussian functions.

First, as an immediate consequence of \cite[Thm.~3.5]{AMP06}, for arbitrary initial data $u_0 \in L^1$, there exists a solution with $u(t,\cdot) \in L^1$ for all $t>0$ via the existence of a positive contraction semigroup in $L^1$. However, the scheme we consider here is most naturally suited for producing approximations in $L^2$, where the question of boundedness of $\phi_t^{\fp}$ arises. In particular, for the choices of polynomial $F$ of main interest, $\phi_t^{\fp} \colon L^2 \to L^2$ can in general be unbounded for any $t>0$ when considered on the unweighted space $L^2$ \cite{LLMR22}.

This difficulty can be overcome by considering the evolution in a suitably weighted space $L^2(\rho)$. Sufficient conditions on $\rho$, depending on $F$, that ensure the boundedness of $\phi_t^{\fp}: L^2(\rho)\to L^2(\rho)$ are given in \cite{GMM2017}. These conditions concern conservative fields $F=-\nabla V$ with confining potentials $V$, for which solutions converge to a stationary state. Our setting is more general, as $F$ need not be conservative or admit a stationary state. The conditions in \cite{GMM2017} would be unnecessarily restrictive for our purposes. We therefore adopt a more tailored approach, exploiting the fact that our method considers initial data $u_0$ given by linear combinations of Gaussian functions, or approximated by such combinations. We focus on assumptions sufficient to handle this class of initial data.

\begin{theorem}\label{thm:semigroupboundedness}
  Let $b \geq 0$ such that
  \begin{equation}\label{eq:bcond}
    \sup_{x \in \R^d} \bigl\{ \kappa  b^2 \abs{x}^2 + b x \cdot F(x) - \textstyle\frac12\displaystyle \nabla \cdot F(x)  \bigr\} \leq \omega_b < \infty.
  \end{equation}
  Then for $\rho_b$ as defined in \eqref{eq:rhoa}, we have
  \begin{equation*}
    \norm{\phi_t^{\fp} u_0}_{L^2(\rho_b)} \leq e^{\omega_b t} \norm{ u_0 }_{L^2(\rho_b)}.
  \end{equation*}
  The same statement holds in the case $\kappa = 0$ corresponding to the pure advection semigroup.
\end{theorem}

\begin{proof}
  We adapt the argument in \cite[Lemma~3.8]{GMM2017}. Let $v(t, \cdot) = \rho_{b/2} \phi_t^{\fp} u_0$. Then
  $$
    \begin{aligned}
      \frac12 \frac{\mathrm d}{\mathrm dt} \norm{v}_{L^2}^2 & = \bigl\langle \rho_{b/2} \nabla \cdot \bigl(\kappa \nabla (\rho_{b/2}^{-1} v ) - F \rho_{b/2}^{-1} v \bigr), v\bigr\rangle                                                   \\
                                                            & = - \kappa \norm{\nabla v}_{L^2}^2  + \bigl\langle \bigl(\kappa  b^2 \abs{x}^2 + b x \cdot F(x) - \textstyle\frac12\displaystyle \nabla \cdot F(x) \bigr) v, v \bigr\rangle .
    \end{aligned}
  $$
  With \eqref{eq:bcond}, the statement follows by a Gronwall argument, since $\rho_{b/2}^2 = \rho_b$.
\end{proof}

\begin{cor}\label{cor:gaussiansemigroupboundedness}

  Let $u_0 \in V_k=\text{span}\set{g_{w_i}}_{i=1}^k$ where the Gaussian functions are linearly independent, and let $a^*> 0$ be as in \eqref{eq:astar}. If there exists $b\geq 0$ such that $0\leq b < a^*$ and \eqref{eq:bcond} holds, then
  \begin{equation*}
    \norm{\phi_t^{\fp} u_0}_{L^2} \leq C(b,V_k) e^{\omega_b t} \norm{ u_0 }_{L^2},
  \end{equation*}
  with $C(b, V_k)$ as in \eqref{eq:upper-constant-norm-eq}.
\end{cor}

\begin{proof}
  This follows immediately from \Cref{thm:semigroupboundedness} and the fact that since $0\leq b< a^*$, \Cref{lemma:norm-equivalence-guassian-subspaces} guarantees that the $L^2$- and $L^2(\rho_b)$-norms are equivalent on $V_k$.
\end{proof}

\section{Main scheme}\label{sec:main_scheme}
\subsection{Splitting and re-approximations}
\label{sec:lie_projection}

Our time integrator relies on a splitting scheme based on a decomposition of the right-hand side of the Fokker--Planck equation \eqref{eq:FP} into a
diffusion operator with affine drift and the remaining non-affine drift term. To this end, the drift $F$ is decomposed into an affine component
$$
  F_{\degone}(x)=\Gamma x+\beta,
$$
with $\Gamma\in\bR^{d\times d}$ and $\beta\in\bR^d$, and a generally non-affine remainder
$$
  F_{\deghigher}=F-F_{\degone}.
$$
With this splitting, we define the operators $\mathcal A$ and $\mathcal B$ by
\begin{equation*}
  \mathcal A v=\kappa\Delta v - \nabla\cdot(F_{\degone}v), \qquad \mathcal Bv=-\nabla\cdot(F_{\deghigher}v).
\end{equation*}
The numerical scheme now alternates between the flows $\phi_t^{\df}$ and $\phi_t^{\adv}$ associated to
\begin{equation}\label{eq:nonlinear_transport}
  \partial_t u = \mathcal A u, \qquad
  \partial_t u = \mathcal B u,
\end{equation}
respectively, and we denote by $\Phi_t^{\df}$ and $\Phi_t^{\adv}$ the numerical approximation to $\phi_t^{\df}$ and $\phi_t^{\adv}$.

\begin{figure}[tbp]
  \begin{center}
    \centering
    \begin{tikzpicture}[
        scale=0.94,
        transform shape,
        x=1cm,y=1cm,
        >=Latex,
        gaussian/.style={draw=green!45!black, thick},
        diffarrow/.style={->, thick, black!75},
        chararrow/.style={->, very thick, blue!70!black},
        projectarrow/.style={->, very thick, orange!85!black, dashed},
        note/.style={font=\footnotesize, align=center},
        glabel/.style={font=\footnotesize, align=center, text=green!45!black}
      ]
      \begin{scope}[shift={(0,0)}]
        \draw[gaussian] (0,0) ellipse (0.55 and 0.28);
        \draw[gaussian] (0,0) ellipse (0.32 and 0.16);
        \draw[gaussian, rotate=18] (0.75,0.48) ellipse (0.45 and 0.22);
        \draw[gaussian, rotate=-14] (0.55,-0.45) ellipse (0.38 and 0.18);
        \node[glabel] at (0.35,-1.18) {$u^n\in\bG_{k_n}$};
      \end{scope}

      \draw[diffarrow] (1.35,0) -- (3.15,0)
      node[midway, above, note] {$\phi_\tau^{\df}$};

      \begin{scope}[shift={(4.20,0.02)}]
        \draw[gaussian] (0,0) ellipse (0.68 and 0.34);
        \draw[gaussian] (0,0) ellipse (0.39 and 0.20);
        \draw[gaussian, rotate=18] (0.78,0.48) ellipse (0.55 and 0.27);
        \draw[gaussian, rotate=-14] (0.58,-0.45) ellipse (0.47 and 0.22);
        \node[glabel] at (0.35,-1.18) {$\phi_\tau^{\df}(u^n)\in\bG_{k_n}$};
      \end{scope}

      \draw[diffarrow] (5.65,0) -- (7.85,0)
      node[midway, above, note] {$\mathcal R_{\delta_n}$};

      \begin{scope}[shift={(8.80,0)}]
        \draw[gaussian] (0,0) ellipse (0.72 and 0.36);
        \draw[gaussian] (0,0) ellipse (0.42 and 0.21);
        \draw[gaussian, rotate=18] (0.60,0.42) ellipse (0.58 and 0.29);
        \node[glabel] at (0.10,-1.18) {$\mathcal R_{\delta_n}\circ\phi_\tau^{\df}(u^n)$};
      \end{scope}

      \draw[diffarrow] (8.05,-0.20) .. controls (6.10,-2.10) and (2.80,-1.05) .. (1.65,-2.80)
      node[pos=0.50, below, note] {$\Phi_{\tau,n}^{\adv}$};

      \fill[blue!70!black] (9.58,0.59) circle (1.5pt);
      \node[note, blue!70!black] at (9.14,0.58) {$x(0)$};
      \draw[chararrow] (9.58,0.59) .. controls (12.00,-0.45) and (12.00,-2.45) .. (9.42,-3.32);
      \fill[blue!70!black] (9.42,-3.32) circle (1.5pt);
      \node[note, blue!70!black] at (9.00,-3.18) {$x(\tau)$};
      \node[note, blue!70!black] at (12.50,-1.50) {$\dot x=F_{\deghigher}(x)$};

      \draw[->, thick, black!75, dashed] (8.90,-1.55) -- (8.90,-2.80)
      node[midway, left, note] {$\phi_\tau^{\adv}$};

      \begin{scope}[shift={(8.80,-3.70)}]
        \draw[red!70!black, thick]
        plot[smooth cycle, tension=0.75] coordinates {
            (-0.55,0.08) (-0.35,0.55) (0.20,0.72) (0.78,0.42)
            (0.62,-0.05) (0.88,-0.50) (0.15,-0.68) (-0.48,-0.42)
          };
        \draw[red!70!black, thick]
        plot[smooth cycle, tension=0.75] coordinates {
            (-0.28,0.02) (-0.15,0.30) (0.20,0.39) (0.52,0.23)
            (0.42,-0.03) (0.52,-0.28) (0.12,-0.38) (-0.25,-0.22)
          };
        \node[note, red!70!black] at (0.10,-1.10)
        {$\phi_\tau^{\adv}\circ\mathcal R_{\delta_n}\circ\phi_\tau^{\df}(u^n)$\\generally outside the LCG ansatz};
      \end{scope}

      \draw[projectarrow] (7.95,-3.70) -- (2.40,-3.70)
      node[midway, above, note] {mass-constrained $L^2$ fitting};

      \begin{scope}[shift={(1.50,-3.70)}]
        \draw[gaussian] (-0.05,0.04) ellipse (0.40 and 0.20);
        \draw[gaussian, rotate=-20] (-0.45,0.42) ellipse (0.34 and 0.17);
        \draw[gaussian, rotate=25] (0.62,0.48) ellipse (0.38 and 0.18);
        \draw[gaussian, rotate=16] (-0.50,-0.38) ellipse (0.32 and 0.15);
        \draw[gaussian, rotate=-15] (0.52,-0.45) ellipse (0.36 and 0.17);
        \node[glabel] at (0.15,-1.10) {$u^{n+1}=\Phi_{\tau,n}^{\adv}\circ\mathcal R_{\delta_n}\circ\phi_\tau^{\df}(u^n)$\\$\in\bG_{k_{n+1}}$};
      \end{scope}
    \end{tikzpicture}
  \end{center}

  \caption{One affine-first Lie--Trotter step, $u^{n+1}=\Phi_{\tau,n}^{\adv}\circ\mathcal R_{\delta_n}\circ\phi_\tau^{\df}(u^n)$. Green ellipses represent Gaussian components. Given an input $g\in \bG_{k_n}$,  the exact affine Fokker--Planck propagation $\phi_\tau^{\df}(g)$ gives an element in $\bG_{k_n}$, thus preserving the number of components. Recompression $\cR_{\delta_n}$ is applied immediately afterward, potentially decreasing that number. The solid curved arrow denotes our numerical transport map $\Phi_{\tau,n}^{\adv}$. The blue curve illustrates a characteristic of $F_{\deghigher}$. The dashed black arrow leads to the exact transport target shown in red. The orange arrow represents mass-constrained $L^2$ fitting of a Taylor approximation of the transported density.}\label{fig:scheme}
\end{figure}
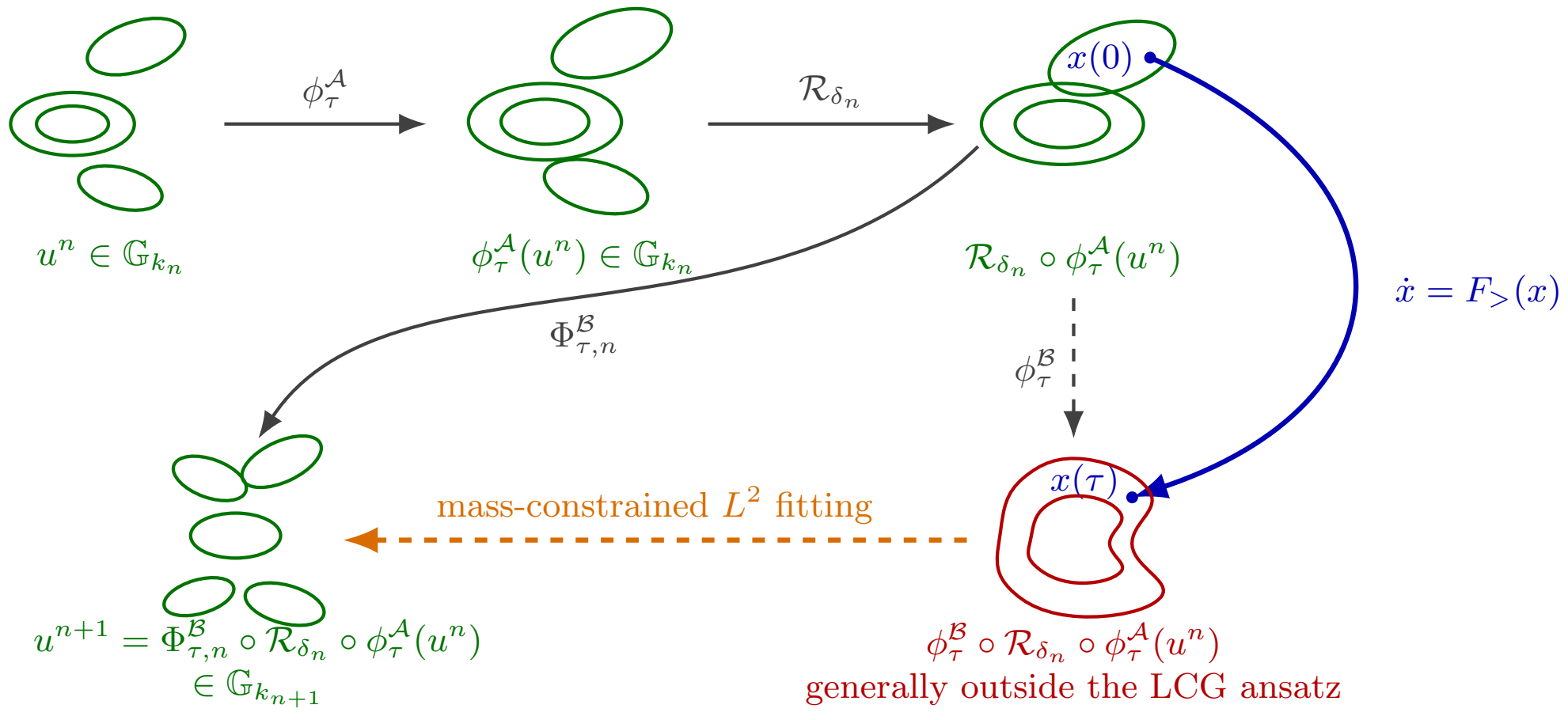

As shown below in \Cref{lem:affine_gaussian_preservation}, we have $\phi_t^{\df}(\bG_k) \subseteq \bG_k$ for each $k$, and the action of $\phi_t^{\df}$ can be evaluated in closed form; we thus use
$$
  \Phi_t^{\df} = \phi_t^{\df}.
$$
In contrast, for any $t>0$ and $k \in \N$, the action of $\phi_t^{\adv}$ generally leaves $\bG_k$ for such initial data. Thus suitable approximations $\Phi_t^{\adv}$ of the transport semi-group are required such that $\Phi_t^{\adv}(g) \in \bG_{\bar k_g(t)}$ for some $\bar k_g(t)$ depending on the given $g \in \bG_k$. Assuming for a moment that such an approximation is available, we can perform time integration based on standard splitting schemes. In the following, we focus on the classical Lie--Trotter splitting for simplicity. Setting a time step $\tau >0$, we define the corresponding one-step map as
\[
  \Phi_{\tau,n} = \Phi^{\adv}_{\tau,n}\circ\mathcal R_{\delta_n}\circ \phi^{\df}_{\tau},
\]
where $\Phi^{\adv}_{\tau,n}$ denotes the action of our numerical transport flow solver at step $n$ for a time step of size $\tau$, and $\mathcal R_{\delta_n}$ is what we call a recompression operator with tolerance $\delta_n$. Our approximation of the exact solution at time $t_n = n\tau$ is given by
\begin{equation}\label{eq:lie_scheme}
  u^{n+1} = \Phi_{\tau,n}(u^n), \quad u^0=u_0, \quad \forall n\in \bN_0.
\end{equation}
We next summarize the main ideas underlying the construction of $\cR_\delta$ and $\Phi_{\tau, n}^{\adv}$ (see Section~\ref{sec:adaptive}). \Cref{fig:scheme} summarizes one time step of our approach.

The role of the recompression operator is to eliminate redundancies in the LCG approximation and to avoid an unnecessarily large number of terms. Given a tolerance $\delta>0$ and an input LCG $u_{c,w} \in \bG_{k}$ with $k$ Gaussians, $\mathcal R_{\delta}(u_{c,w}) \in \bG_{\tilde k_\delta}$ provides an approximation of $u_{c,w}$ with $\tilde k_\delta \leq k$ such that $\norm{\mathcal R_{\delta}(u_{c,w}) - u_{c,w}}_{L^2} \leq \delta$ subject to the constraint that $\mathcal R_{\delta}(u_{c,w})$ has unit mass. Finding the LCG that realizes the minimal possible $\tilde k_\delta$ is generally infeasible, since it amounts to finding the global minimizers of a sequence of high-dimensional non-convex optimization problems. For particular objectives and error criteria, a variety of methods has been proposed for reduction of Gaussian mixtures; see, for example, \cite{CrouseWillettPattipatiSvensson2011,BeylkinMonzonYang2019}.
We present in Section~\ref{sec:adaptive_recompression} a greedy strategy adapted to our purposes that provides a reduction in practice and guarantees that the prescibed tolerance $\delta$ is reached.

In the remainder of this section, we provide in \Cref{sec:affine} the details of the exact integrator $\phi_t^{\df}(u_0)$ when $u_0 \in \bG_k$. In \Cref{sec:transport} we explain our specific construction of $\Phi_t^{\adv}$. We anticipate here that it is essentially based on three building blocks:
\begin{enumerate}
  \item \textbf{Taylor approximation of the transported density:} The exact semi-group $\phi_t^{\adv}$ of \eqref{eq:nonlinear_transport} can be expressed through the characteristic flow map generated by $F_{\deghigher}$. Since the flow map is generally not available in closed form, we approximate the transported density by a Taylor expansion in time whose fitting quantities can be evaluated in closed form.
  \item \textbf{Projection onto $\mathbb{G}_k$.} For a given $k\in \mathbb N$, we approximate the Taylor-expanded transported density by an element of $\mathbb{G}_k$ through mass-constrained $L^2$ fitting.
  \item \textbf{Adaptive increase in the number of Gaussians:} We evaluate the $L^2$-error in fitting the Taylor approximation of the transported density. We increase the number of Gaussians $k$ if this error exceeds a certain target accuracy.
\end{enumerate}

\subsection{Exact affine Fokker--Planck integrator}\label{sec:affine}
In this section, we show that $\phi_t^{\df}$ maps LCGs to LCGs while preserving the mass of each Gaussian component. This result is folklore in the literature, but its proof is not usually explicitly provided. We thus recall the proof for convenience and self-completeness. Note that when the affine drift vanishes, $\phi_t^{\df}$ reduces to the heat flow.
\begin{lemma}\label[lemma]{lem:affine_gaussian_preservation}
  For every $t\geq 0$, the solution to
  \begin{align}
    \label{eq:FP-linear-drift}
     & \partial_t u = \mathcal A u, \qquad u(0,\cdot)=\sum_{i=1}^k c_i(0) g_{w_i(0)}  \;\in \bG_k,
  \end{align}
  with $k\in \bN$, $\mathcal A u = \kappa\Delta u - \nabla\cdot\bigl(F_{\degone}u\bigr)$ and $F_{\degone}=\Gamma x+\beta$, is
  $$
    \phi_t^{\df}(u_0) = \sum_{i=1}^k c_i(t) g_{w_i(t)} \;\in \bG_k
  $$
  with
  $$
    c_i(t)=c_i(0)\left(\frac{\det P_i(t)}{\det P_i(0)}\right)^{1/2},
    \qquad
    m_i(t)=e^{t\Gamma}m_i(0)+\int_0^t e^{(t-s)\Gamma}\beta\,\d s,
  $$
  and $U_i(t)$ defined for $i=1,\ldots,k$ by the unique Cholesky factorization $P_i(t)=U_i(t)^\top U_i(t)$ of
  $$
    P_i(t)=\left(e^{t\Gamma}P_i^{-1}(0)e^{t\Gamma^\top}+2\kappa\int_0^t e^{(t-s)\Gamma}e^{(t-s)\Gamma^\top}\,\d s\right)^{-1}.
  $$
\end{lemma}
\begin{proof}
  By linearity of \eqref{eq:FP-linear-drift}, it suffices to consider a single Gaussian initial datum $u(0,\cdot) = c_0 g_{m_0, U_0} \in \bG_1$ with $P_0 = U_0^\top U_0$.

  The spatial Fourier transform of $u(t,\cdot) \in L^1\cap L^2$ is given by
  $$
    \hat u(t, \xi)=(2\pi)^{-d/2}\int_{\R^d}u(t,x)e^{-i\xi^\top x}\,\d x, \xi \in \bR^d.
  $$
  The Fourier transform is also well defined for our initial datum in $\bG_1\subset L^1\cap L^2$, and
  $$
    \hat u(0,\xi)
    =
    c_0\, \widehat{g_{m_0,U_0}}(\xi)
    =
    c_0\,\det(P_0)^{-1/2}
    \exp\left(-i\xi^\top m_0-\frac12\xi^\top P_0^{-1}\xi\right),
    \quad  \xi \in \bR^d.
  $$
  Applying the Fourier transform to the equation gives
  $$
    \partial_t\hat u-(\Gamma^\top\xi)\cdot\nabla_{\xi}\hat u
    +i\xi^\top\beta\,\hat u+\kappa|\xi|^2\hat u=0.
  $$
  For every $\xi_0\in \bR^d$, the associated characteristic satisfies
  $$
    \frac{\d \xi(s)}{\d s}=-\Gamma^\top\xi(s),
    \qquad \xi(0)=\xi_0,
  $$
  hence $\xi(s)=e^{-s\Gamma^\top}\xi_0$. Along such a characteristic,
  $$
    \frac{\d \hat u(s,\xi(s))}{\d s}
    =
    -\left(i\xi(s)^\top\beta+\kappa|\xi(s)|^2\right)\hat u(s,\xi(s)).
  $$
  Therefore
  \begin{align*}
    \hat u(t,\xi)
     & =
    \hat u(0,e^{t\Gamma^\top}\xi)
    \exp\left(
    -i\xi^\top\int_0^t e^{(t-s)\Gamma}\beta\,\d s
    -\kappa\xi^\top\int_0^t e^{(t-s)\Gamma}e^{(t-s)\Gamma^\top}\,\d s\,\xi
    \right) \\
     & =
    c_0\,\det(P_0)^{-1/2}
    \exp\left(
    -i\xi^\top m(t)-\frac12\xi^\top P(t)^{-1}\xi
    \right).
  \end{align*}
  On the other hand,
  $$
    \widehat{c(t)g_{m(t),U(t)}}(\xi)
    =
    c(t)\det(P(t))^{-1/2}
    \exp\left(
    -i\xi^\top m(t)-\frac12\xi^\top P(t)^{-1}\xi
    \right).
  $$
  Noting that $c(t)\det(P(t))^{-1/2}=c_0\det(P_0)^{-1/2}$ completes the proof.
\end{proof}

\subsection{Transport integrator by Taylor expansion}\label{sec:transport}
Let
$$
  v_0(x)=\sum_{j=1}^{k_0}c_{0,j}g_{w_{0,j}}(x)
$$
be the datum at the start of a transport substep. Let $\psi_t$ denote the flow generated by the nonlinear drift (see the blue curve in \Cref{fig:scheme}),
\begin{equation}\label{eq:charac}
  \frac{d}{dt}\psi_t(x)=F_{\deghigher}(\psi_t(x)),\qquad \psi_0(x)=x.
\end{equation}
Assume that $F_{\deghigher}$ is polynomial and that there exists $t_0>0$ such that, for every $t\in[0,t_0]$, its flow $\psi_t$ is defined on $\mathbb R^d$ and is a $C^1$ diffeomorphism onto $\Omega_t=\psi_t(\mathbb R^d)$, with $J_t=\det \D\psi_t>0$. We denote its inverse on $\Omega_t$ by $\psi_t^{-1}$. Setting $u_{\adv}(t) = e^{t\mathcal B} v_0$, we obtain
\begin{equation}\label{eq:u_adv}
  u_{\adv}(t,\psi_t(x))=v_0(x)
  \exp\left(-\int_0^t\nabla\cdot F_{\deghigher}(\psi_s(x))\,\d s\right).
\end{equation}
Equivalently, on $\Omega_t$,
\begin{equation}\label{eq:nonlinear_transport_conservative}
  u_{\adv}(t,y)=v_0(\psi_{t}^{-1}(y))\det \D\left(\psi_{t}^{-1}\right)(y),
\end{equation}
with $u_{\adv}$ extended by zero outside $\Omega_t$.

Direct evaluation of \eqref{eq:nonlinear_transport_conservative} requires the inverse flow map, which is generally unavailable in closed form for nonlinear polynomial drifts. We therefore fit the transported density with an LCG $u_{c,w}=\sum_{i=1}^k c_i g_{w_i}\in\mathbb{G}_k$ in the $L^2$ sense. Other choices of the underlying Hilbert space, such as $H^1$, can be treated in the same manner. The objective can be decomposed as
\begin{equation*}
  \norm{u_{\adv}(t)-u_{c,w}}^2
  =I_{c,w}
  -2 I_{\mathrm{\times}}
  +  I_{\adv},
\end{equation*}
where
\begin{equation*}
  I_{c,w} = \norm{u_{c,w}}^2 , \quad
  I_{\mathrm{\times}} =  \langle u_{\adv}(t),u_{c,w}\rangle, \quad
  I_{\adv} = \norm{u_{\adv}(t)}^2,
\end{equation*}
and where we write $\norm{\cdot} = \norm{\cdot}_{L^2}$ and $\langle\cdot,\cdot\rangle=\langle\cdot,\cdot\rangle_{L^2}$.

The first term is explicit. Indeed,
$$
  I_{c,w}=\sum_{i,j}c_ic_j\int_{\R^d}g_{w_i}(x)g_{w_j}(x)\,\d x=c^\top G(w)c,
$$
where $G_{ij}(w)=\langle g_{w_i},g_{w_j}\rangle$. Writing $P_i=U_i^\top U_i$ and $P_j=U_j^\top U_j$, we have
\begin{equation*}
  G_{ij}(w)=\exp\left(-\frac12\left(m_i^\top P_i m_i+m_j^\top P_j m_j-m_c^\top P_c m_c\right)\right)\frac{(2\pi)^{d/2}}{\sqrt{\det P_c}},
\end{equation*}
with $P_c=P_i+P_j$, $m_c=P_c^{-1}(P_i m_i+P_j m_j)$.

It remains to approximate $I_{\mathrm{\times}}$ and $I_{\adv}$. For $r\in\mathbb N_+$, define the Taylor approximation of the transported density by
\begin{equation}\label{eq:transport_taylor_target}
  \widetilde u_{\adv}^{(r)}(t)
  =\sum_{\ell=0}^r\frac{t^\ell}{\ell!}\mathcal B^\ell v_0.
\end{equation}
Provided that $u_{\adv}$ is sufficiently regular in time, we have $\mathcal B^\ell v_0=\left.\partial_t^\ell u_{\adv}(t)\right|_{t=0}, \ell=0,\ldots,r$. Thus, $\widetilde u_{\adv}^{(r)}$ is the Taylor polynomial of the transported density in time.
Since $F_{\deghigher}$ is polynomial, each term is a finite sum of polynomials times Gaussians. We use this same approximation in both the cross term $I_{\mathrm{\times}}$ and the squared norm $I_{\adv}$.
The next two lemmas give approximations of these terms and bounds for the corresponding errors.
\begin{lemma}[Approximation of $I_{\mathrm{\times}}$]\label[lemma]{lem:transport_cross}
  Let $r\in \mathbb N_+$ and fix the Gaussian parameters $(c,w)$. Setting
  \begin{equation}\label{eq:b_vector_transport}
    b_i(w;t)=\langle g_{w_i},\widetilde u_{\adv}^{(r)}(t)\rangle
    =\sum_{\ell=0}^r\frac{t^\ell}{\ell!}
    \langle g_{w_i},\mathcal B^\ell v_0\rangle,
  \end{equation}
  we have
  \begin{equation*}
    \widehat I_{\mathrm{\times}}(t)                           =c^\top b(w;t)
    =\langle\widetilde u_{\adv}^{(r)}(t),u_{c,w}\rangle,                       \qquad
    |I_{\mathrm{\times}}(t)-\widehat I_{\mathrm{\times}}(t)|  \lesssim t^{r+1},
  \end{equation*}
  for $0<t\leq t_0$.
\end{lemma}

\begin{lemma}[Approximation of $I_{\adv}$]\label[lemma]{lem:transport_norm}
  Let $r\in \mathbb N_+$ and assume $u_{\adv}\in C^{r+1}([0,t_0];L^2)$. Define
  \begin{equation}\label{eq:hat_I_nl}
    \widehat I_{\adv}(t)=\norm{\widetilde u_{\adv}^{(r)}(t)}^2=\sum_{\ell,m=0}^r\frac{t^{\ell+m}}{\ell!m!}\langle\mathcal B^\ell v_0,\mathcal B^m v_0\rangle.
  \end{equation}
  Moreover, for $0<t\leq t_0$,
  \begin{align*}
    \norm{u_{\adv}(t)-\widetilde u_{\adv}^{(r)}(t)}
                                       & \leq\frac{t^{r+1}}{(r+1)!}
    \sup_{0\leq s\leq t_0}\norm{\partial_s^{r+1}u_{\adv}(s)},       \\
    |I_{\adv}(t)-\widehat I_{\adv}(t)| & \lesssim t^{r+1}.
  \end{align*}
\end{lemma}

The proofs of Lemmas \ref{lem:transport_cross} and \ref{lem:transport_norm} are given in Appendix~\ref{sec:transport_taylor_appendix}, where we provide the polynomial recurrences and explicit cross coefficients. Since each $\mathcal B^\ell v_0$ is a finite sum of polynomials times Gaussians, all coefficients in the expansions of $\widehat I_{\mathrm{\times}}(t)$ and $\widehat I_{\adv}(t)$ in powers of $t$ can be evaluated using moments of the form
\begin{equation}\label{eq:poly_gaussian_moment}
  \int_{\R^d}q(x)g_{w_i}(x)g_{w_j}(x)\,\d x,
\end{equation}
where $q$ is polynomial. Since a product of Gaussians is a constant multiple of a Gaussian, these moments can be evaluated exactly using Gaussian moment formulas.

For fixed Gaussian parameters $(c,w)$, Lemmas~\ref{lem:transport_cross} and \ref{lem:transport_norm}  yield
$$
  \norm{\widetilde u_{\adv}^{(r)}(\tau)-u_{c,w}}_{L^2}^2
  =c^\top G(w)c-2c^\top b(w;\tau)+\widehat I_{\adv}(\tau),
$$
which approximates $\norm{u_{\adv}(\tau)-u_{c,w}}_{L^2}^2$ with an error of order $\mathcal O(\tau^{r+1})$.

In what follows, we fix the Taylor order $r\geq1$ and the step size $\tau$, and abbreviate $b(w;\tau)$ as $b(w)$. We define the squared fitting error with respect to $\widetilde u_{\adv}^{(r)}(\tau)$ by
\begin{equation*}
  \mathcal E(c,w)
  =\norm{u_{c,w}-\widetilde u_{\adv}^{(r)}(\tau)}_{L^2}^2
  =c^\top G(w)c-2c^\top b(w)+\widehat I_{\adv}(\tau).
\end{equation*}
Since $\widehat I_{\adv}(\tau)=\norm{\widetilde u_{\adv}^{(r)}(\tau)}_{L^2}^2$ is independent of $(c,w)$, minimizing $\mathcal E$ is equivalent to minimizing
\begin{equation}\label{eq:obj_fun}
  \mathcal J(c,w)=c^\top G(w)c-2c^\top b(w).
\end{equation}
We therefore determine the Gaussian parameters by solving the mass-constrained optimization problem
\begin{equation}\label{eq:opt_problem}
  \inf_{\substack{c\in\R^k,\;w\in\mathbb W^k\\\mathrm{a}(w)^\top c=1}} \mathcal J(c,w),
\end{equation}
where $\mathrm{a}(w)$ collects the masses of the Gaussian components:
\begin{equation*}
  \mathrm{a}_i(w)=\int_{\R^d}g_{w_i}(x)\,\dx
  =\frac{(2\pi)^{d/2}}{\sqrt{\det P_i}},
  \qquad
  \mathrm{a}(w)=(\mathrm{a}_1(w),\ldots,\mathrm{a}_k(w))^\top.
\end{equation*}
The output of the approximate transport integrator $\Phi_\tau^{\adv}$ is defined as the LCG corresponding to the solution of \eqref{eq:opt_problem} obtained by the method described in Section~\ref{sec:adaptive}.

\section{Optimization and Adaptivity}\label{sec:adaptive}
The goal of this section is to construct an LCG approximation of $\widetilde u_{\adv}^{(r)}(\tau)$ with an appropriate number of components by solving \eqref{eq:opt_problem}. Specifically, for a prescribed fitting error tolerance $\varepsilon_{\text{tol}}$, we seek an LCG $\Phi_\tau^{\adv}v_0$ satisfying
$$
  \norm{\Phi_\tau^{\adv}v_0-\widetilde u_{\adv}^{(r)}(\tau)}_{L^2}\leq \varepsilon_{\text{tol}}.
$$
To begin with, we describe how the Gaussian parameters are determined for a fixed number of components in the following subsection.

\subsection{Gaussian fitting with a fixed number of components}\label{sec:adaptive_optimization}
The coefficients $c$ enter the objective function \eqref{eq:obj_fun} quadratically and can therefore be eliminated before optimizing the Gaussian shape parameters.

\subsubsection{Elimination of linear coefficients and analytic derivatives} For fixed Gaussian shape parameters $w$ with linearly independent Gaussian components, so that $G(w)$ is positive definite, the coefficients $c$ in \eqref{eq:opt_problem} are determined as a function of $w$ by
\begin{equation*}
  c(w) = \argmin_{\mathrm{a}(w)^\top c=1} c^\top G(w)c-2c^\top b(w).
\end{equation*}
The first-order optimality conditions with multiplier $\lambda(w)$ read $G(w)c(w)=b(w)-\lambda(w) \mathrm{a}(w)$ and $\mathrm{a}(w)^\top c(w)=1$, with the explicit solution
\begin{equation}\label{eq:mass_weights}
  c(w)=G(w)^{-1}\bigl(b(w)-\lambda(w) \mathrm{a}(w)\bigr),\qquad
  \lambda(w)=\frac{\mathrm{a}(w)^\top G(w)^{-1}b(w)-1}{\mathrm{a}(w)^\top G(w)^{-1}\mathrm{a}(w)}.
\end{equation}
Substituting these expressions into \eqref{eq:obj_fun} gives the reduced optimization problem
\begin{equation}\label{eq:reduced_problem}
  \inf_{w\in\mathbb W^k}\mathcal J(w), \quad  \mathcal J(w)=\mathcal J(c(w),w)=-b(w)^\top c(w)-\lambda(w).
\end{equation}
The corresponding squared fitting error indicator is
\begin{equation}\label{eq:squared_fitting_error_reduced}
  \mathcal E (w) = \mathcal{E} (c(w), w) = \mathcal J(w)+\widehat I_{\adv}(\tau).
\end{equation}

We next derive explicit expressions for the derivatives of $\mathcal J(w)$. Using the coefficient optimality conditions, differentiation with respect to the shape parameters gives the Fr\'echet derivative
\begin{equation*}
  \D \mathcal J
  =c^\top(\D G)c-2c^\top(\D b)
  +2\lambda(\D\mathrm{a})^\top c.
\end{equation*}
Here, $c=c(w)$ and $\lambda=\lambda(w)$ are given by \eqref{eq:mass_weights}.

Since $U$ is upper triangular, its free entries are $U_{ij}$ with $1\leq i\leq j\leq d$. The derivatives of the Gaussian components are
\[
  \partial_{m_i} g_{m,U}
  =\bigl(U^\top U(x-m)\bigr)_i g_{m,U},\qquad
  \partial_{U_{ij}} g_{m,U}
  =-\bigl(U(x-m)\bigr)_i(x_j-m_j)g_{m,U},\quad i\leq j.
\]
Thus the entries of $\D G$, $\D b$ and $\D\mathrm{a}$ can be expressed as integrals of products of polynomials with Gaussians that can be evaluated explicitly using the Gaussian moments in \eqref{eq:poly_gaussian_moment}.

\subsubsection{Levenberg--Marquardt iteration}
We solve the nonconvex problem \eqref{eq:reduced_problem} locally using a variable-projection Levenberg--Marquardt method. At each evaluation of the reduced objective, the coefficients $c(w)$ are determined by \eqref{eq:mass_weights}.

At the current shape parameters $w$, let $h = (h_1,\ldots,h_k)$ be a trial increment, where $h_i$ corresponds to component $i$. By linearizing the fitted LCG, we obtain
$u_{c(w+h),w+h}
  \approx u_{c(w),w}
  +\D_w u_{c(w),w}\,h$.
Substituting this linearization into the expression for $\mathcal J(w+h)$ yields the Gauss--Newton model
\[
  \mathcal J(w+h)
  \approx\mathcal J(w)+\D \mathcal J(w) h
  + h^\top H h,
\]
where $H=(\D_w u_{c(w),w})^* (\D_w u_{c(w),w})$ with the Hilbert adjoint of $\D_w u_{c(w),w} \colon \R^{\dim \bW^k} \to L^2$; note that the entries of $H$ again reduce to integrals of products of polynomials and Gaussians.

To control the size of the update, we introduce a scaled damping penalty \cite{more1978}. Let
$B_i
  =(\D_{w_{i}}g_{w_i})^*(\D_{w_{i}}g_{w_i})$
denote the derivative Gramian blocks, and consider the block-diagonal damping matrix
\begin{equation*}
  D=\operatorname{diag}(c_1^2B_1,\ldots,c_k^2B_k).
\end{equation*}
With the coefficients held fixed, the associated quadratic penalty satisfies
\begin{align*}
  h^\top D h
  =\sum_{i=1}^k c_i^2 h_i^\top B_i h_i
  =\sum_{i=1}^k\bignorm{c_i (\D_{w_i} g_{w_i} )h_{i}}_{L^2}^2.
\end{align*}
Thus, the penalty measures the sum of the squared linearized changes in the individual components, accounting for the different sensitivities and coupling of the parameters within each Gaussian.

\begin{remark}
  In practice, small or vanishing coefficients can make $D$ nearly singular. To retain damping for the corresponding components, we replace each factor $c_i^2$ in $D$ by
  \begin{equation*}
    \max\left\{c_i^2,\;10^{-3}
    \frac{\max_j c_j^2G_{jj}}{G_{ii}}\right\}.
  \end{equation*}
  Since $c_i^2G_{ii}=\norm{c_i g_{w_i}}_{L^2}^2$, the replacement sets a lower bound of \(10^{-3}\max_j\|c_jg_{w_j}\|_{L^2}^2\) on the squared component sizes used for damping. In what follows, $D$ denotes this modified damping matrix. The inner products defining $H$ and $D$ are evaluated using analytic Gaussian moments, with the coefficients determined by \eqref{eq:mass_weights}.
\end{remark}

The trial step minimizes the damped local model
\begin{equation}\label{eq:lm_local_model}
  \min_{h}
  \left\{\mathcal J(w)+\D\mathcal J(w) h + h^\top H h+ \mu h^\top D h\right\},
  \qquad \mu>0.
\end{equation}
Its first-order optimality condition is
\begin{equation*}
  (H+\mu D) h =-\tfrac12\D\mathcal J(w)^\top,\qquad \mu>0,
\end{equation*}
which we solve after block Cholesky scaling by $D$.

A feasible trial step is accepted if it decreases the objective, in which case we update $w \gets w+h$. Trials with singular precision factors or Gram matrices are rejected. The damping parameter $\mu$ is initialized to $10^{-3}$ and adjusted using the ratio of the actual objective reduction to the reduction predicted by the Gauss--Newton model:
\[
  \rho = - \frac{\mathcal J(w)-\mathcal J(w+h)}{\D\mathcal J(w) h + h^\top H h}.
\]
For an accepted step, we increase $\mu$ when $\rho<\frac12$, decrease it when $\rho>\frac12$, and leave it unchanged when $\rho=\frac12$. For a rejected step, we increase $\mu$.

To reduce accumulated fitting error, we use an inner tolerance $0<\varepsilon_{\mathrm{in}}<\varepsilon_{\mathrm{tol}}$. The iteration terminates when
\[
  \mathcal E(w)\leq\varepsilon_{\mathrm{in}}^2,
\]
or when the maximum number of iterations is reached.

\subsection{Adaptive refinement}\label{sec:adaptive_refinement}
Let $w^\ast$ denote the shape parameters returned by the Levenberg--Marquardt iteration. Using the fitting error indicator \eqref{eq:squared_fitting_error_reduced}, the fit is accepted if
\begin{equation}\label{eq:stop_condition}
  \mathcal E(w^\ast)\leq\varepsilon_{\mathrm{tol}}^2.
\end{equation}
Since the transport step generally maps a Gaussian to a non-Gaussian function, the prescribed tolerance may not be attainable with the current number of components. In this case, we enrich the approximation by splitting a selected Gaussian component, indexed by $i^\ast$, into two child components. We then reoptimize all shape parameters and coefficients using the procedure described in Section~\ref{sec:adaptive_optimization}. Optimization and enrichment alternate until the stopping criterion \eqref{eq:stop_condition} is satisfied or the maximum number of splits $s_{\max}$ is reached. Algorithm~\ref{alg:adaptive_transport} summarizes this greedy-type adaptive procedure.
\begin{algorithm}[htp]
  \caption{Adaptive transport fitting}\label{alg:adaptive_transport}
  \begin{algorithmic}[1]
    \Require $v_0$, step size $\tau$, tolerance $\varepsilon_{\mathrm{tol}}$, maximum number of splits $s_{\max}$
    \Ensure Best LCG approximation $v_{\mathrm{best}}$
    \State Construct $\widetilde u_{\adv}^{(r)}(\tau)$ and initialize the shape parameters $w$
    \State Set $v_{\mathrm{best}}=u_{c(w),w}$ and $\mathcal E_{\mathrm{best}}=\mathcal E(w)$
    \For{$j=0,\ldots,s_{\max}$}
    \State Optimize $w$ by Levenberg--Marquardt iteration with tolerance $\varepsilon_{\text{in}}$ and coefficients from \eqref{eq:mass_weights}
    \State Evaluate the squared fitting error $\mathcal E(w)$ with updated $w$
    \If{$\mathcal E(w)<\mathcal E_{\mathrm{best}}$}
    \State Set $v_{\mathrm{best}}=u_{c(w),w}$ and $\mathcal E_{\mathrm{best}}=\mathcal E(w)$
    \EndIf
    \If{$\mathcal E_{\mathrm{best}}\leq\varepsilon_{\mathrm{tol}}^2$ or $j=s_{\max}$}
    \State \Return $v_{\mathrm{best}}$
    \EndIf
    \State Select component $i_\ast$ for splitting, enrich the Gaussian dictionary, and update $w$.
    \EndFor
  \end{algorithmic}
\end{algorithm}

We next describe how to select the component to split in Section~\ref{sec:component_splitting} and how to construct the enriched Gaussian dictionary in Section~\ref{sec:split_direction}. The enriched dictionary along with its updated coefficients provides the initial guess for the subsequent optimization. In the remainder of this subsection, $w$ and $c=c(w)$ denote the parameters of the current fit, and we abbreviate $g_i=g_{w_i}$, with $w_i=(m_i,U_i)$. 

\subsubsection{Component splitting indicator}\label{sec:component_splitting}
When the fit fails, the transport has deformed some components into shapes that no single Gaussian can follow.
To describe what can be captured by Gaussians, for any $w = (m,U)\in \bW$, we define the tangent space
\[
  \mathcal T_{w} = \mathcal T_{m,U}
  =
  \Span\left\{
  g_{m,U},
  \partial_{m_\alpha}g_{m,U},
  \partial_{U_{\beta\gamma}}g_{m,U}
  \colon
  1\leq\alpha\leq d, 1\leq\beta\leq\gamma\leq d
  \right\},
\]
and set $\mathcal T_i = \mathcal T_{w_i}$.

One explicit Euler step of the transport maps a component to
\[
  c_ig_i+\tau\adv(c_ig_i)
  =\mathrm e^{\tau\adv}(c_ig_i)+\mathcal O(\tau^2).
\]
A single Gaussian can follow this perturbation only by changing its coefficient, mean, and covariance, that is, only within the corresponding tangent space. For each component, the indicator measures the part of $\tau\adv(c_ig_i)$ outside this tangent space.

The following characterization is the real-valued counterpart of \cite[Lemma~3.1]{LasserLubich2020}. For completeness, we give its proof in Appendix~\ref{sec:proof_tangent}.
\begin{lemma}\label{lem:tangent}
  Let $\mathcal P_{\leq2}$ denote the space of real polynomials on $\R^d$ of total degree at most two. For every $(m,U)\in\mathbb W$,
  \[
    \mathcal T_{m,U}=g_{m,U}\mathcal P_{\leq2},
    \qquad \dim\mathcal T_{m,U}=\frac{(d+1)(d+2)}2.
  \]
\end{lemma}

Let $\Pi_i$ denote the $L^2(\R^d)$-orthogonal projection onto $\mathcal T_i$. We define the normal deformation and the \emph{splitting indicator}
\begin{equation}\label{eq:eta_indicator}
  \Psi_i=(\id-\Pi_i)\,\adv g_i,
  \qquad
  \eta_i=\tau\,\norm{c_i\Psi_i}_{L^2}=\tau\,|c_i|\,\norm{\Psi_i}_{L^2}.
\end{equation}
With $P_i=U_i^\top U_i$, we have $\adv g_i=g_i\,q_i$ with the polynomial
\[
  q_i=F_{\deghigher}^\top P_i(x-m_i)-\nabla\cdot F_{\deghigher}
\]
of degree at most $\deg F_{\deghigher}+1$. We evaluate $\eta_i$ for all components and split a component with index
\[
  i^\ast
  \in
  \argmax_{1\leq i\leq k}\eta_i .
\]

The indicator measures the best single-Gaussian fit of the transported component. Let
\[
  E_i(\tau)=\min_{|(c,w)-(c_i,w_i)|\le M\tau}\norm{\mathrm e^{\tau\adv}(c_ig_i)-c\,g_w}_{L^2}
\]
for a fixed $M>0$ large enough. Since a Gaussian near $c_ig_i$ can only follow the part of $\tau\adv(c_ig_i)$ in $\mathcal T_i$, linearization gives $E_i(\tau)=\eta_i+\mathcal O(\tau^2)$. By linearity of $\mathrm e^{\tau\adv}$,
\[
  \inf_{u\in\bG_k}\norm{\mathrm e^{\tau\adv}u_{c,w}-u}_{L^2}
  \le\sum_{i=1}^kE_i(\tau)
  =\sum_{i=1}^k\eta_i+\mathcal O(\tau^2).
\]
The indicator is a selection criterion, not an estimator of the fitting error, since components may cancel.

\begin{remark}\label{rem:indicator_pullback}
  If $\deg F_{\deghigher}\le3$, then $\nabla\cdot F_{\deghigher}\in\mathcal P_{\le2}$, so the divergence part of $\adv g_i$ lies in $\mathcal T_i$ and is removed by $\id-\Pi_i$. By \eqref{eq:generic_cross_taylor} in the proof of Lemma~\ref{lem:transport_cross},
  \[
    g_i(\psi_\tau(x))
    =
    g_i(x)
    \bigl(
    p_{\times,i}^{(r)}(x;\tau)
    +\mathcal O(\tau^{r+1})
    \bigr),
  \]
  and the first-order term of $g_i\,p_{\times,i}^{(r)}$ is $\tau F_{\deghigher}\cdot\nabla g_i$. It follows that
  \[
    \eta_i=|c_i|\,\bigl\|(\id-\Pi_i)\bigl(g_i\,p_{\times,i}^{(r)}\bigr)\bigr\|_{L^2}+\mathcal O(\tau^2).
  \]
  The right-hand side uses only the polynomial $p_{\times,i}^{(r)}$ already computed for the cross terms, and it is the version used in all experiments of Section~\ref{sec:numerical_examples}.
\end{remark}

It remains to explain how to compute $\eta_i$ and its variant in Remark~\ref{rem:indicator_pullback}. Both are of the form $\norm{(\id-\Pi_i)(g_ih)}_{L^2}$ for a polynomial $h$, namely $h=q_i$ or $h=p_{\times,i}^{(r)}$. By Lemma~\ref{lem:tangent}, $g_i[h]_{\le2}\in\mathcal T_i$, where $[\cdot]_{\le2}$ and $[\cdot]_{\ge3}$ denote the parts of total degree at most two and at least three, so that $(\id-\Pi_i)(g_ih)=(\id-\Pi_i)(g_i[h]_{\ge3})$. The degree decomposition is not orthogonal after multiplication by $g_i$, for example, in one dimension,
\[
  \left\langle x^3g_i,xg_i\right\rangle_{L^2}=\int_{\R}x^4g_i(x)^2\,\dx>0,
\]
so the projection is still needed; it is easy to see that it can be evaluated explicitly in terms of integrals of products of Gaussians with polynomials.

\subsubsection{Splitting direction}\label{sec:split_direction}
We now construct the two child components of the selected component $i^\ast$. Set $\Sigma_i=(U_i^\top U_i)^{-1}$ and fix a splitting parameter $\gamma\in(0,1)$. Given a nonzero direction $v\in\R^d$, define the displacement
\[
  \zeta
  =
  v\sqrt{\frac{\gamma}{v^\top\Sigma_{i^\ast}^{-1}v}}.
\]
The child components have means $m_{i^\ast}\pm\zeta$ and share the covariance
\[
  \Sigma_{\mathrm c}
  =
  \Sigma_{i^\ast}-\zeta\zeta^\top.
\]
Since $\zeta^\top\Sigma_{i^\ast}^{-1}\zeta=\gamma<1$, the matrix $\Sigma_{\mathrm c}$ is positive definite. The coefficients are chosen so that each child carries half the mass of the parent component. As a consequence, the combined children preserve the mass, mean, and covariance of the parent.

To select $v$, we compare coordinate splits in the standardized variables $y=U_{i^\ast}(x-m_{i^\ast})$. For $v=U_{i^\ast}^{-1}e_j$, the child Gaussians have means $\pm\sqrt\gamma\,e_j$ and covariance $I_d-\gamma e_je_j^\top$. For the resulting deviation in $L^2$-norm, we have the following asymptotic expansion.

\begin{prop}\label{prop:symmetric_split_expansion}
  For $y \in \R^d$ and $\gamma \in (0,1)$, let $\bar P_\gamma = (I_d - \gamma e_j e_j^\top)^{-\frac12} $ and
  \[
    \bar \sigma_\gamma(y) = \frac{1}{2\sqrt{1-\gamma}} \biggl[ \exp\Bigl( - \frac12 \bigabs{ \bar P_\gamma (y + \sqrt{\gamma} e_j)}^2   \Bigr)  + \exp\Bigl( - \frac12 \bigabs{\bar P_\gamma  (y - \sqrt{\gamma} e_j)}^2   \Bigr)    \biggr].
  \]
  Then there exists $C>0$ depending only on $d$ such that for all $\gamma \in (0, \frac12)$,
  \begin{equation}\label{eq:symmetric_split_expansion}
    \int_{\R^d} \biggabs{ \bar \sigma_\gamma(y)
      -  \exp\left( - \frac12\abs{y}^2 \right) \left[
        1-\frac{\gamma^2}{12}(y_j^4-6y_j^2+3)
        \right] }^2 \,dy \leq C \gamma^6 .
  \end{equation}
\end{prop}

\begin{proof}
  This follows from Mehler's formula \cite[Chapter 6]{AAR}, which for the probabilist's Hermite polynomials $\mathrm{He}_n$, $n\in\N_0$, yields
  \[
    \frac{\bar \sigma_\gamma(y)}{\exp\left( - \frac12\abs{y}^2 \right) } = \sum_{n=0}^\infty \frac{\gamma^n \mathrm{He}_{2n}(1)}{(2n)!} \mathrm{He}_{2n}(y_j) = 1 - \frac{\gamma^2}{12} \mathrm{He}_{4}(y_j) + \frac{\gamma^3}{45} \mathrm{He}_6(y_j) + \ldots,
  \]
  where $\mathrm{He}_{4}(y_j) = y_j^4 - 6 y_j^2 + 3$; another application of Mehler's formula yields that the resulting remainder can be bounded in $L^2$ to give the right-hand side of \eqref{eq:symmetric_split_expansion}.
\end{proof}

This expansion motivates a direction score based on the overlap between $\Psi_{i^\ast}$ and the leading perturbation induced by the split. Since $\Psi_{i^\ast}$ is orthogonal to $g_{i^\ast}q$ for all $q\in\mathcal P_{\le2}$, only the quartic term contributes to this overlap. We therefore define the scores
\[
  S_j
  =
  \left|
  \int_{\R^d}
  \Psi_{i^\ast}(x)g_{i^\ast}(x)
  \bigl[U_{i^\ast}(x-m_{i^\ast})\bigr]_j^4
  \,\dx
  \right|,
  \qquad j=1,\ldots,d.
\]
As an additional heuristic restriction, we retain only indices $j$ such that
$F_{\deghigher,j}\not\equiv0$ or
$\partial_{x_j}F_{\deghigher,\ell}\not\equiv0$ for some $\ell$.
Among the retained indices, we select $j^\ast$ with the largest score $S_{j^\ast}$ and set
\[
  v=U_{i^\ast}^{-1}e_{j^\ast}.
\]
If all retained scores vanish numerically, we instead take $v=e_j$, where $j$ is a retained index maximizing $(\Sigma_{i^\ast})_{jj}$. If all component indicators vanish numerically, we select a component maximizing $|c_i|\operatorname{tr}(\Sigma_i)$ and choose $v$ along an eigenvector corresponding to its largest covariance eigenvalue.

Finally, the two child components are expressed in the parametrization \eqref{eq:lcg} using a Cholesky factorization of $\Sigma_{\mathrm c}^{-1}$. The enriched approximation serves as the initial guess for the subsequent optimization. Its accuracy may therefore influence the convergence of the Levenberg--Marquardt iteration, while the accuracy of the accepted approximation is controlled by the prescribed tolerance $\varepsilon_{\text{tol}}$.

\subsection{Recompression}\label{sec:adaptive_recompression}
The affine Fokker--Planck update preserves the number of Gaussian components. Repeated transport enrichment may nevertheless produce redundant components. We therefore recompress LCGs by retaining a subset of the Gaussian components and refitting their coefficients up to a prescribed $L^2$ tolerance. To this end, we use an orthogonal least squares approach \cite{ChenBillingsLuo1989,ChenCowanGrant1991}.
For an input
$$
  u_{c,w}=\sum_{i=1}^k c_i g_{w_i},
  \qquad c\in\R^k,\quad w=(w_1,\dots,w_k)\in\mathbb W^k,
$$
with $\mathrm{a}(w)^\top c=1$ (where $\mathrm{a}(w)$ is the mass vector defined in Section~\ref{sec:transport}) and a recompression tolerance $\delta\geq0$, the recompression operator
$\mathcal R_\delta$ produces $u_{\widehat c,\widehat w}=\mathcal R_\delta(u_{c,w})\in\bG_{\widehat k}$, where $1\leq\widehat k\leq k$, such that
$$
  \norm{u_{c,w}-u_{\widehat c,\widehat w}}_{L^2}\leq\delta,
  \qquad
  \int_{\R^d}u_{\widehat c,\widehat w}(x)\,\d x=1.
$$
To reduce the computational cost of transport, we apply recompression at each time step, after the affine Fokker--Planck flow and before the transport step in \eqref{eq:lie_scheme}.

Recompression selects Gaussian shapes from $w$ and refits their coefficients by mass-constrained orthogonal least squares. Define the normalized components by
$$
  \tilde g_i=\frac{g_{w_i}}{\norm{g_{w_i}}_{L^2}}.
$$
For each nonempty support $S\subseteq\{1,\ldots,k\}$, let
$$
  \tilde{\mathrm{a}}(S)=\left(\int_{\R^d}\tilde g_i(x)\,\d x\right)_{i\in S}
  =\left(\frac{\mathrm{a}_i(w)}{\norm{g_{w_i}}_{L^2}}\right)_{i\in S},
$$
and define the best approximation error subject to a unit mass constraint by
\begin{equation}\label{eq:recompressionerror}
  E(S)= \min \biggl\{  \Bignorm{u_{c,w}-\sum_{i\in S}\widehat c_i^S \tilde g_i }_{L^2} \colon  \widehat c^S \in \R^S , \tilde{\mathrm{a}}(S)^\top \widehat c^S  = 1 \biggr\}.
\end{equation}
Adding components enlarges the feasible set and hence cannot increase $E(S)$. Starting with $S=\emptyset$, the forward pass repeatedly chooses
$$
  j_*\in\operatorname*{arg\,min}_{j\notin S}E(S\cup\{j\}),
  \qquad S\leftarrow S\cup\{j_*\},
$$
until $E(S)\leq\delta$ or all components have been selected.

We extend this to a forward-backward greedy selection following \cite{Zhang2011}. Once the tolerance is met, a backward pass tests every one-component deletion. Among deletions whose refitted error remains at most $\delta$, it accepts one with the smallest error. This is repeated until no admissible deletion remains or only one component is left.
In exact arithmetic, the full support satisfies $E(\{1,\ldots,k\})=0$ because the input has unit mass. If the numerical projection on the full support fails to meet the tolerance, the input is returned unchanged.

To evaluate the constrained projections in \eqref{eq:recompressionerror}, form the normalized Gram matrix
$$
  \tilde{G}_{ij}=\langle\tilde g_i,\tilde g_j\rangle_{L^2},\qquad
  \tilde{G}=N^{-1}G(w)N^{-1},\qquad N=\operatorname{diag}(\norm{g_{w_1}}_{L^2},\ldots,\norm{g_{w_k}}_{L^2}).
$$
Writing $\tilde{G}^{\frac12}$ for its positive semidefinite square root, set
$$
  X=(\tilde{G}^{\frac12})_{:,S},\qquad y=\tilde{G}^{\frac12}Nc.
$$
For $a_0=\tilde{\mathrm{a}}(S)/\|\tilde{\mathrm{a}}(S)\|_2^2$, and with the columns of $Z$ forming an orthonormal basis of $\ker\tilde{\mathrm{a}}(S)^\top$, a minimizing coefficient vector is given by
$$
  \widehat c^S=a_0+Z(XZ)^+(y-X a_0),
$$
where the pseudoinverse is computed by singular value decomposition. The residual norm
$\|y-X\widehat c^S\|_2$ equals $E(S)$. For the accepted support $S$, the coefficients in the original Gaussian parametrization are
$$
  \widehat w=(w_i)_{i\in S},\qquad
  \widehat c=\left(\frac{\widehat c_i^S}{\norm{g_{w_i}}_{L^2}}\right)_{i\in S},\qquad
  \widehat k=|S|,
$$
so that $u_{\widehat c,\widehat w}=\sum_{i\in S}\widehat c_i^S\tilde g_i$.

\section{Error analysis}\label{sec:error_analysis}
In this section, we estimate the $L^2$ error of the Lie--Trotter scheme \eqref{eq:lie_scheme}, separating splitting, Taylor truncation, fitting, and recompression errors.

With $\phi_t^{\fp}$ the exact flow of \eqref{eq:FP} and $\phi_t^{\df}=\mathrm e^{t\mathcal A}$, $\phi_t^{\adv}=\mathrm e^{t\mathcal B}$, the exact and numerical Lie--Trotter maps are
\begin{align}
  \phi_\tau^{\mathrm{lie}} & = \phi_{\tau}^{\adv}\circ\phi_{\tau}^{\df}=\mathrm e^{\tau\mathcal B}\mathrm e^{\tau\mathcal A},\notag \\
  \Phi_{\tau,n}            & = \Phi_{\tau,n}^{\adv}\circ\mathcal R_{\delta_n}\circ\phi_{\tau}^{\df}.\label{eq:full_numerical_map}
\end{align}
The component maps in \eqref{eq:full_numerical_map} are defined in Sections~\ref{sec:main_scheme} and~\ref{sec:adaptive}. When recompression is omitted, $\mathcal R_{\delta_n}$ is the identity. Since the split generators are unbounded, the analysis requires control of commutators in weighted Sobolev spaces \cite{jahnke_lubich2000}.

Define $X_t=L^2(J_t^{-1})$, where $J_t(x)=\det \D\psi_t(x)$ is the Jacobian determinant of the characteristic flow $\psi_t$ (see \eqref{eq:charac} for the definition) and $J_t^{-1}(x)=1/J_t(x)$. Since $J_0=1$, we have $X_0=L^2$. The semigroup property $\psi_{\theta+t}=\psi_\theta\circ\psi_t$ and the chain rule give
$$
  J_{\theta+t}(y)=J_\theta(\psi_t(y))J_t(y).
$$
Changing variables $x=\psi_t(y)$ in \eqref{eq:nonlinear_transport_conservative} therefore yields, for any $\theta,t\geq0$ with $\theta+t\leq T$ and $f\in X_{\theta+t}$,
\begin{align}\label{eq:isometry_identity}
  \|\mathrm e^{t\mathcal B}f\|_{X_\theta}^2
  =\int_{\R^d}\frac{|f(y)|^2}{J_\theta(\psi_t(y))J_t(y)}\,\d y
  =\norm{f}_{X_{\theta+t}}^2.
\end{align}
Hence $\mathrm e^{t\mathcal B}$ is an isometry from $X_{\theta+t}$ to $X_\theta$. In particular, it is an isometry from $X_t$ to $L^2$ when $\theta=0$. For our error analysis, we use this identity together with the following stability assumption for the affine Fokker–Planck semigroup.

\begin{assumption}\label{assmp:weighted_stability}
  {\rm(i)} $F:\R^d\to\R^d$ is a polynomial vector field with $\deg F\le q$ for some integer $q\ge1$. \\
  {\rm(ii)} For every $0\leq t\leq T$, the characteristic flow $\psi_t$ is a $C^1$ diffeomorphism from $\R^d$ onto its image, with $J_t>0$. \\
  {\rm(iii)} There exist $\tau_0>0$ and constants $m_0,m_2\geq0$ such that, for every $0\leq \tau \leq\tau_0$ and $0\leq \theta\leq T$,
  \begin{align}
    \norm{\mathrm e^{\tau \mathcal A}v}_{X_\theta}
     & \leq \mathrm e^{m_0\tau}\norm{v}_{X_\theta},\label{eq:affine_l2_stability}                                \\
    \norm{\mathrm e^{\tau \mathcal A}v}_{H^2(\langle x\rangle^{2q}J_\theta^{-1})}
     & \leq \mathrm e^{m_2\tau}\norm{v}_{H^2(\langle x\rangle^{2q}J_\theta^{-1})},\label{eq:affine_h2_stability}
  \end{align}
  for every $v$ in the respective space, where $\langle x\rangle=(1+|x|^2)^{1/2}$.
\end{assumption}

\begin{example}[Double-well and radial drifts]\label{ex:weighted_stability}
  We verify \eqref{eq:affine_l2_stability} and \eqref{eq:affine_h2_stability} for two drifts used in our numerical tests, both with $q=3$. Set
  $$
    \rho_{t,0}=J_t^{-1},\qquad\rho_{t,2}=\langle x\rangle^{2q}J_t^{-1}.
  $$
  A sufficient condition for the required affine Fokker--Planck bounds is
  $$
    C_\ell=\sup_{0\leq t\leq T}
    \left\|
    \frac{\kappa\Delta\rho_{t,\ell}
    +F_{\degone }\cdot\nabla\rho_{t,\ell}
    -(\nabla\cdot F_{\degone })\rho_{t,\ell}}
    {\rho_{t,\ell}}
    \right\|_{L^\infty}<\infty,
    \qquad \ell\in\{0,2\}.
  $$
  For fixed $t$, the weighted energy argument in the proof of Theorem~\ref{thm:semigroupboundedness} (see also \cite[proof of Lemma~3.8]{GMM2017}) yields
  \begin{align*}
    \frac{\d}{\d \tau}\norm{\mathrm e^{\tau \mathcal A}v}_{L^2(\rho_{t,\ell})}^2
     & =-2\kappa\int_{\R^d}\rho_{t,\ell}|\nabla(\mathrm e^{\tau \mathcal A}v)|^2\,\dx        \\
     & \quad+\int_{\R^d}\bigl(\kappa\Delta\rho_{t,\ell}+F_{\degone }\cdot\nabla\rho_{t,\ell}
    -(\nabla\cdot F_{\degone })\rho_{t,\ell}\bigr)|\mathrm e^{\tau \mathcal A}v|^2\,\dx      \\
     & \leq C_\ell\norm{\mathrm e^{\tau \mathcal A}v}_{L^2(\rho_{t,\ell})}^2.
  \end{align*}
  By Gronwall's inequality, \eqref{eq:affine_l2_stability} holds with $m_0=C_0/2$. Since $F_{\degone }(x)=\Gamma x+\beta$, the commutator $[\partial^\alpha,\mathcal A]$ is a constant linear combination of derivatives of order $|\alpha|$ for $1\leq|\alpha|\leq2$. Applying the same energy argument at the $H^2$ level establishes \eqref{eq:affine_h2_stability} with $m_2=(C_2+C_\Gamma)/2$, where $C_\Gamma\geq0$ depends only on $\Gamma$ and $d$.

  \begin{enumerate}[(i)]
    \item For the double-well drift $F(x)=(a x_1(1-x_1^2),-b x_2,\ldots,-b x_d)$, with $a,b>0$, choose
          $$
            F_{\degone }(x)=(a x_1,-b x_2,\ldots,-b x_d),\qquad
            F_{\deghigher }(x)=(-a x_1^3,0,\ldots,0).
          $$
          Here $\Gamma=\operatorname{diag}(a,-b,\ldots,-b)$. Thus,
          \begin{align*}
            \psi_t(x)_1 & =\frac{x_1}{\sqrt{1+2atx_1^2}},
            \qquad \psi_t(x)_i=x_i\quad (i\geq2),         \\
            J_t^{-1}(x) & =(1+2atx_1^2)^{3/2}.
          \end{align*}

    \item For the radial drift $F(x)=4(1-|x|^2)x$, take
          $$
            F_{\degone }(x)=4x,\qquad F_{\deghigher }(x)=-4|x|^2x.
          $$
          Here $\Gamma=4I_d$. Thus,
          $$
            \psi_t(x)=\frac{x}{\sqrt{1+8t|x|^2}},
            \qquad J_t^{-1}(x)=
            \left(1+8t|x|^2\right)^{(d+2)/2}.
          $$
  \end{enumerate}
  To check $C_0,C_2<\infty$, use $\Delta\rho/\rho=\Delta\log\rho+|\nabla\log\rho|^2$ and, for $c\geq0$,
  \begin{align*}
    \nabla\log(1+c|x|^2)
     & =\frac{2cx}{1+c|x|^2},                                \\
    \Delta\log(1+c|x|^2)
     & =\frac{2cd}{1+c|x|^2}-\frac{4c^2|x|^2}{(1+c|x|^2)^2}.
  \end{align*}
  The corresponding one-dimensional identities apply to the factor $1+2atx_1^2$. With $c\in\{1,2at,8t\}$ and $0\leq c|x|^2/(1+c|x|^2)\leq1$, these formulas show that $\nabla\log\rho_{t,\ell}$, $\Delta\log\rho_{t,\ell}$, and $F_{\degone }\cdot\nabla\log\rho_{t,\ell}$ are bounded uniformly in $x$ and $0\leq t\leq T$. Thus $C_0,C_2<\infty$ for both drifts, which establishes the stability required in \eqref{eq:affine_l2_stability} and \eqref{eq:affine_h2_stability}.
\end{example}
\begin{remark}
  If $J_\theta^{-1}(x)\le C_T\rho_b(x)$ uniformly for $0\le\theta\le T$, then Theorem~\ref{thm:semigroupboundedness} gives
  $$
    \norm{\phi_t^{\fp}f}_{X_\theta}\le C_T^{1/2}e^{\omega_b t}\norm{f}_{L^2(\rho_b)}
  $$
  for $f\in L^2(\rho_b)$. The splitting analysis below additionally requires weighted $H^2$ bounds.
\end{remark}
The next result isolates the defect caused by the noncommutativity of the transport and affine Fokker--Planck flows.
\begin{prop}[Lie--Trotter splitting error]\label{prop:lie_splitting_error}
  Suppose Assumption~\ref{assmp:weighted_stability} holds. Let $\theta\geq0$ and $0\leq\tau\leq\tau_0$ satisfy $\theta+\tau\leq T$. Assume that the exact solution $v(s)=\phi_s^{\fp}v$ belongs to $H^2(\langle x\rangle^{2q}J_{\theta+\tau-s}^{-1})$ for $0\leq s\leq\tau$, with norms bounded uniformly in $s$. Then
  $$
    \norm{\phi_\tau^{\mathrm{fp}}v-\phi_\tau^{\mathrm{lie}}v}_{X_\theta}
    \le C_{\mathrm{lie}}\tau^2
    \sup_{0\leq s\leq\tau}
    \norm{v(s)}_{H^2(\langle x\rangle^{2q}J_{\theta+\tau-s}^{-1})},
  $$
  where $C_{\mathrm{lie}}$ depends only on $\kappa$, the coefficients of $F_{\degone }$ and $F_{\deghigher }$, $m_0$, $m_2$, and $\tau_0$.
\end{prop}
\begin{proof}
  For our ordering of the two flows, the Duhamel commutator argument in \cite[Theorem~2.1]{iserles_kropielnicka2024} leads to
  $$
    \phi_\tau^{\mathrm{lie}}v-\phi_\tau^{\mathrm{fp}}v=
    \int_0^\tau\int_0^{\tau-s}
    \mathrm e^{(\tau-s)\mathcal B}\mathrm e^{(\tau-s-\sigma)\mathcal A}[\mathcal B, \mathcal A]\mathrm e^{\sigma\mathcal A}\mathrm e^{s(\mathcal A+\mathcal B)}v\,\d \sigma\, \d s,
  $$
  where $[X,Y]=XY-YX$. To obtain this identity, differentiate $\mathrm e^{(\tau-s)\mathcal B}\mathrm e^{(\tau-s)\mathcal A}\mathrm e^{s(\mathcal A+\mathcal B)}v$ with respect to $s$ and use that for $h>0$,
  \[
    [\mathcal B,\mathrm e^{h\mathcal A}]
    =\int_0^h\mathrm e^{(h-\sigma)\mathcal A}[\mathcal B,\mathcal A]\mathrm e^{\sigma\mathcal A}\,\d\sigma.
  \]
  The calculation is first performed for smooth data and extended by approximation using the bounds below.
  The operator $[\mathcal B,\mathcal A]$ has differential order at most $2$, and its coefficients grow at most like $\langle x\rangle^q$. Therefore, for every positive weight $\rho$,
  \begin{align*}
    \norm{[\mathcal B,\mathcal A]v}_{L^2(\rho)}
     & \le C_1\norm{v}_{H^2(\langle x\rangle^{2q}\rho)},
  \end{align*}
  where $C_1>0$ depends only on $\kappa$ and the coefficients of $F_{\degone }$ and $F_{\deghigher }$.

  Using \eqref{eq:isometry_identity}, the two stability estimates, and the commutator bound, we obtain
  \begin{align*}
    \norm{\mathrm e^{(\tau-s)\mathcal B}\mathrm e^{(\tau-s-\sigma)\mathcal A}[\mathcal B,\mathcal A]\mathrm e^{\sigma\mathcal A}\mathrm e^{s(\mathcal A+\mathcal B)}v}_{X_\theta} & =\norm{\mathrm e^{(\tau-s-\sigma)\mathcal A}[\mathcal B,\mathcal A]\mathrm e^{\sigma\mathcal A}\mathrm e^{s(\mathcal A+\mathcal B)}v}_{X_{\theta+\tau-s}} \\
                                                                                                                                                                                  & \le C_1\mathrm e^{m_0(\tau-s-\sigma)}\mathrm e^{m_2\sigma}\norm{v(s)}_{H^2(\langle x\rangle^{2q}J_{\theta+\tau-s}^{-1})}.
  \end{align*}
  With $m=\max\{m_0,m_2\}$, integration over the triangle $0\leq s\leq\tau$ and $0\leq\sigma\leq\tau-s$ yields
  \begin{align*}
    \norm{\phi_\tau^{\mathrm{lie}}v-\phi_\tau^{\mathrm{fp}}v}_{X_\theta}
     & \le \frac{C_1\tau^2}{2}\mathrm e^{m\tau_0}
    \sup_{0\leq s\leq\tau}
    \norm{v(s)}_{H^2(\langle x\rangle^{2q}J_{\theta+\tau-s}^{-1})}.
  \end{align*}
  Thus one can take $C_{\mathrm{lie}}=\frac{C_1}{2}\mathrm e^{m\tau_0}$.
\end{proof}

\begin{remark}
  The analysis is restricted to the Lie--Trotter splitting in which the affine Fokker--Planck step is applied first. Higher-order splitting methods may be treated by analogous commutator arguments, provided suitable bounds hold for the resulting higher-order commutators and for their propagation in the corresponding weighted Sobolev spaces. We do not pursue these extensions here.
\end{remark}

For the numerical part of the analysis, set
$$
  z^n=\mathcal R_{\delta_n}\circ\phi_\tau^{\df}(u^n).
$$
The order-$r$ Taylor approximation \eqref{eq:transport_taylor_target} with local datum $v_0=z^n$ is
$$
  \widetilde u_{\adv,n}^{(r)}(t)
  =\sum_{\ell=0}^r\frac{t^\ell}{\ell!}\mathcal B^\ell z^n.
$$
We consider a completed run with accepted Gaussian fits and assume that there are nonnegative numbers $\varepsilon_n$, $T_n$, and $\delta_n$ such that, for every $0\leq \theta\leq T-t_{n+1}$,
\begin{equation}\label{eq:weighted_error_budgets}
  \begin{aligned}
    \norm{\Phi_{\tau,n}^{\adv}(z^n)-\widetilde u_{\adv,n}^{(r)}(\tau)}_{X_\theta} \le \varepsilon_n,\quad
    \norm{\widetilde u_{\adv,n}^{(r)}(\tau)-\phi_\tau^{\adv}(z^n)}_{X_\theta}\le T_n, \quad
    \norm{z^n-\phi_\tau^{\df}(u^n)}_{X_{\theta+\tau}} \le \delta_n.
  \end{aligned}
\end{equation}
Here $\varepsilon_n$, $T_n$, and $\delta_n$ bound the fitting, Taylor truncation, and recompression errors, respectively. In $L^2$, the fitting error is $\sqrt{\mathcal E(w^\ast)}$ from \eqref{eq:squared_fitting_error_reduced}, evaluated for the local datum $v_0=z^n$. The corresponding weighted bounds are additional assumptions needed for the global analysis. In particular, the recompression bound in $X_{\theta+\tau}$ is required because recompression is followed by transport.
\begin{remark}\label{rem:weighted_control}
  The fitting and recompression bounds in \eqref{eq:weighted_error_budgets} involve weighted norms, whereas the experiments in Section~\ref{sec:numerical_examples} use unweighted $L^2$ norms. These norms coincide for the nonlinear oscillator (see Section~\ref{sec:num_duffing}) and the spinning top (see Section~\ref{sec:num_spinning_top}), since $\nabla\cdot F_{\deghigher}=0$ implies $J_\theta\equiv1$ and hence $X_\theta=L^2$. For the double-well and radial drifts, the formulas in Example~\ref{ex:weighted_stability} give
  \[
    J_\theta^{-1}(x)\leq\rho(x),
    \qquad 0\leq\theta\leq T,
  \]
  with the time-independent polynomial weights
  \[
    \rho(x)=
    \begin{cases}
      (1+2aTx_1^2)^2,                    & \text{double-well drift}, \\
      (1+8T|x|^2)^{\lceil(d+2)/2\rceil}, & \text{radial drift}.
    \end{cases}
  \]
  Performing fitting and recompression in $L^2(\rho)$ therefore enforces the corresponding weighted bounds with the prescribed tolerances. All required inner products can still be evaluated as Gaussian moments of the form \eqref{eq:poly_gaussian_moment}, so the method extends directly to this weighted setting. For the spatially localized LCGs considered here, the additional polynomial weight is expected to have only a limited effect on practical performance.
\end{remark}
The next lemma combines the transport and recompression budgets into the numerical defect of one Lie--Trotter step.
\begin{lemma}[Numerical integrator error]\label{lem:numerical_integrator_error}
  Suppose that \eqref{eq:weighted_error_budgets} holds. Then, for every $0\leq \theta\leq T-t_{n+1}$,
  \[
    \norm{\Phi_{\tau,n}u^n-\phi_\tau^{\mathrm{lie}}u^n}_{X_\theta}
    \le \delta_n+\varepsilon_n+T_n.
  \]
\end{lemma}
\begin{proof}
  Since $\Phi_{\tau,n}u^n=\Phi_{\tau,n}^{\adv}(z^n)$, we obtain
  \begin{align*}
    \norm{\Phi_{\tau,n}u^n-\phi_\tau^{\mathrm{lie}}u^n}_{X_\theta}
     & \leq\norm{\Phi_{\tau,n}^{\adv}(z^n)-\widetilde u_{\adv,n}^{(r)}(\tau)}_{X_\theta} +\norm{\widetilde u_{\adv,n}^{(r)}(\tau)-\phi_\tau^{\adv}(z^n)}_{X_\theta} +\norm{z^n-\phi_\tau^{\df}(u^n)}_{X_{\theta+\tau}} \\
     & \leq\varepsilon_n+T_n+\delta_n
  \end{align*}
  by the triangle inequality and the isometry identity \eqref{eq:isometry_identity}.
\end{proof}

Combining the splitting and numerical defects gives an $L^2$ estimate for one step starting from the numerical state.
\begin{theorem}[One-step error]\label{thm:one_step_error}
  Let Assumption~\ref{assmp:weighted_stability} and \eqref{eq:weighted_error_budgets} hold. Assume that the exact solution $u_n(s)=\phi_s^{\fp}u^n$ belongs to $H^2(\langle x\rangle^{2q}J_{\tau-s}^{-1})$ for $0\leq s\leq\tau$, with norms bounded uniformly in $s$. Then for $0<\tau\leq\tau_0$ and $t_{n+1}\leq T$,
  $$
    \norm{\Phi_{\tau,n}u^n-\phi_\tau^{\mathrm{fp}}u^n}_{L^2}
    \le
    C_{\mathrm{lie}}\tau^2\sup_{0\leq s\leq\tau}
    \norm{u_n(s)}_{H^2(\langle x\rangle^{2q}J_{\tau-s}^{-1})}
    +\delta_n+\varepsilon_n+T_n.
  $$
\end{theorem}
\begin{proof}
  Apply the triangle inequality, Lemma~\ref{lem:numerical_integrator_error} with $\theta=0$, and Proposition~\ref{prop:lie_splitting_error} with $v=u^n$ and $\theta=0$.
\end{proof}

Using the local error estimates and the stability of the exact Lie--Trotter map, we obtain the following global error bound.
\begin{theorem}[Global error]\label{thm:global_error}
  Let $T>0$, $n\tau\leq T$, and define
  $$
    u^0=u_0,\qquad u^{j+1}=\Phi_{\tau,j}u^j,\qquad j=0,\ldots,n-1.
  $$
  Suppose Assumption~\ref{assmp:weighted_stability} and the weighted error budgets \eqref{eq:weighted_error_budgets} hold at every step, and let $0<\tau\leq\tau_0$. Assume that the exact solution $u(t)=\phi_t^{\fp}u_0$ satisfies $u(t_j+s)\in H^2(\langle x\rangle^{2q}J_{(n-j)\tau-s}^{-1})$ for $0\leq j<n$ and $0\leq s\leq\tau$, with norms bounded uniformly in $j$ and $s$. Set
  $$
    M_T=\max_{0\leq j<n}\sup_{0\leq s\leq\tau}
    \norm{u(t_j+s)}_{H^2(\langle x\rangle^{2q}J_{(n-j)\tau-s}^{-1})}.
  $$
  Then
  \begin{align*}
    \norm{u^n-\phi_{n\tau}^{\mathrm{fp}}u_0}_{L^2}
    \le \mathrm e^{m_0T}
    \left(C_{\mathrm{lie}}T M_T\tau
    +\sum_{j=0}^{n-1}(\delta_j+\varepsilon_j+T_j)\right).
  \end{align*}
\end{theorem}
\begin{proof}
  By repeated application of \eqref{eq:isometry_identity} and the $X_\theta$ stability of the affine Fokker--Planck semigroup,
  \begin{equation}\label{eq:lie_power_stability}
    \norm{(\phi_\tau^{\mathrm{lie}})^k f}_{L^2}
    \leq \mathrm e^{m_0k\tau}\norm{f}_{X_{k\tau}},
    \qquad k\tau\leq T.
  \end{equation}
  A telescoping decomposition yields
  \begin{align*}
    u^n-\phi_{n\tau}^{\mathrm{fp}}u_0
    =\sum_{j=0}^{n-1}(\phi_\tau^{\mathrm{lie}})^{n-1-j}
    \left[
      (\phi_\tau^{\mathrm{lie}}-\phi_\tau^{\mathrm{fp}})u(t_j)
      +(\Phi_{\tau,j}u^j-\phi_\tau^{\mathrm{lie}}u^j)
      \right].
  \end{align*}
  Set $\theta_j=(n-1-j)\tau$. Applying \eqref{eq:lie_power_stability}, Proposition~\ref{prop:lie_splitting_error} with $v=u(t_j)$ and $\theta=\theta_j$, and Lemma~\ref{lem:numerical_integrator_error} with $\theta=\theta_j$ gives
  \begin{align*}
    \norm{u^n-\phi_{n\tau}^{\mathrm{fp}}u_0}_{L^2}
     & \leq\sum_{j=0}^{n-1}\mathrm e^{m_0\theta_j}
    \left(
    C_{\mathrm{lie}}\tau^2
    \sup_{0\leq s\leq\tau}
    \norm{u(t_j+s)}_{H^2(\langle x\rangle^{2q}J_{\theta_j+\tau-s}^{-1})}
    +\delta_j+\varepsilon_j+T_j
    \right)                                        \\
     & \le \mathrm e^{m_0T}
    \left(C_{\mathrm{lie}}T M_T\tau
    +\sum_{j=0}^{n-1}(\delta_j+\varepsilon_j+T_j)\right). \qedhere
  \end{align*}
\end{proof}

\section{Numerical experiments}\label{sec:numerical_examples}
In this section, to assess the performance of the proposed method, we consider four test examples with different polynomial drifts. In each example, we monitor the evolution of the number of Gaussian components $k(t)$, the fitting error in the transport step, the relative $L^2$ error with respect to a reference solution, and mass conservation. Given a reference solution $u_{\mathrm{ref}}$, we define the relative $L^2$ error at time $t_n$ by
\begin{equation*}
  \mathcal E_{L^2}(t_n)=\frac{\norm{u^n-u_{\mathrm{ref}}(t_n)}_{L^2}}
  {\norm{u_{\mathrm{ref}}(t_n)}_{L^2}}.
\end{equation*}
All transport fitting and recompression are performed in the unweighted $L^2$ norm (see Remark~\ref{rem:weighted_control} for the relation to the weighted bounds assumed in Section~\ref{sec:error_analysis}).
All experiments use the Taylor order $r=3$ and the inner fitting tolerance $\varepsilon_{\mathrm{in}}=0.1\varepsilon_{\mathrm{tol}}$. They were implemented in Julia~1.12.6 and performed on a MacBook Pro equipped with a 10-core Apple M4 processor and 24 GB of memory.

For the examples in Sections~\ref{sec:num_double_well} and~\ref{sec:num_duffing}, we also compare the adaptive scheme with runs using a fixed number $k$ of Gaussian components. In these runs, each time step consists of the exact affine step followed by the same mass-constrained Levenberg--Marquardt fit as in the adaptive scheme, but without enrichment or recompression, so that the number of components remains $k$ throughout.

\subsection{Two-dimensional double-well problem}\label{sec:num_double_well}
Double-well potentials provide a simple model of bistable systems, in which noise can induce transitions between two stable states separated by an energy barrier. The associated Fokker--Planck equation describes the redistribution of probability mass between the two wells. In this example, we consider the drift field
\begin{equation*}
  F(x_1,x_2)=\begin{pmatrix}x_1-x_1^3\\-x_2\end{pmatrix},
\end{equation*}
with affine part $F_{\degone}=(x_1,-x_2)^\top$ and nonlinear residual $F_{\deghigher}=(-x_1^3,0)^\top$. We set the diffusion coefficient to $\kappa=0.3$ and take the initial density to be the single Gaussian $u_0(x)=\frac{1}{2\pi}\exp(-\frac12|x|^2)$. Although the problem is separable in the two coordinates, we solve it directly in two dimensions without exploiting this structure, thereby assessing the effectiveness of the proposed method on the full two-dimensional problem. The solution is computed up to $T=1$ with the fitting tolerance $\varepsilon_{\mathrm{tol}}=10\tau^2$ and the squared recompression tolerance $\delta^2 = \frac12\varepsilon_{\mathrm{tol}}^2$. The reference solution is computed by solving the two one-dimensional marginal equations independently on $[-6,6]$, using Chang--Cooper finite-volume fluxes \cite{CHANG19701} on a uniform grid of $1001$ points, zero-flux boundary conditions, and Crank--Nicolson time stepping with time step size $2.5\times10^{-5}$. The two-dimensional reference density is then obtained as the product of the marginal references.

\begin{figure}[tbp]
  \centering
  \includegraphics[width=0.98\textwidth]{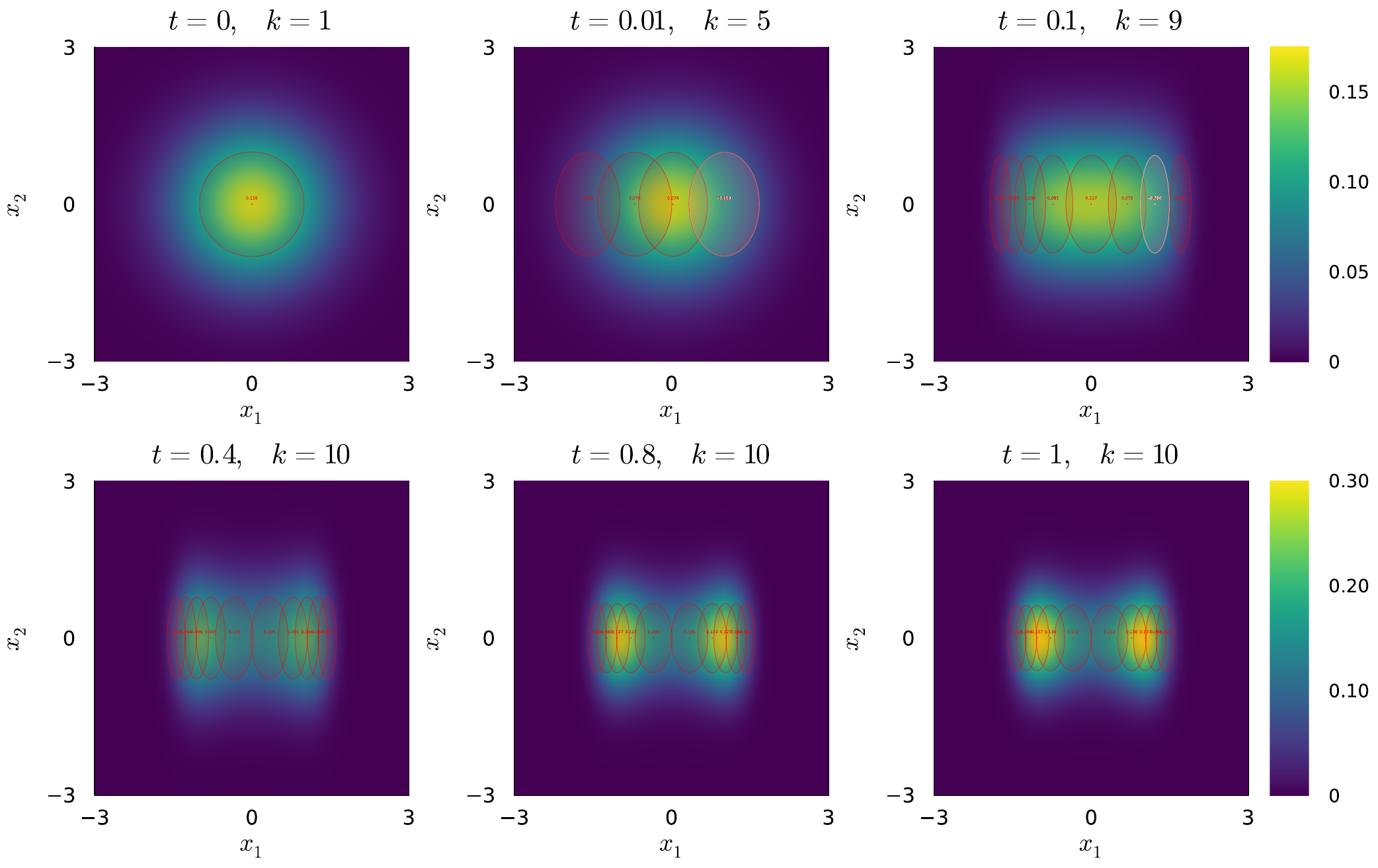}
  \caption{Snapshots of the LCG solution for the double-well problem with $\tau=10^{-3}$. Red and white ellipses represent Gaussian components with positive and negative coefficients, respectively.}\label{fig:dw_snapshots}
\end{figure}
\begin{figure}[tbp]
  \centering
  \includegraphics[width=0.98\textwidth]{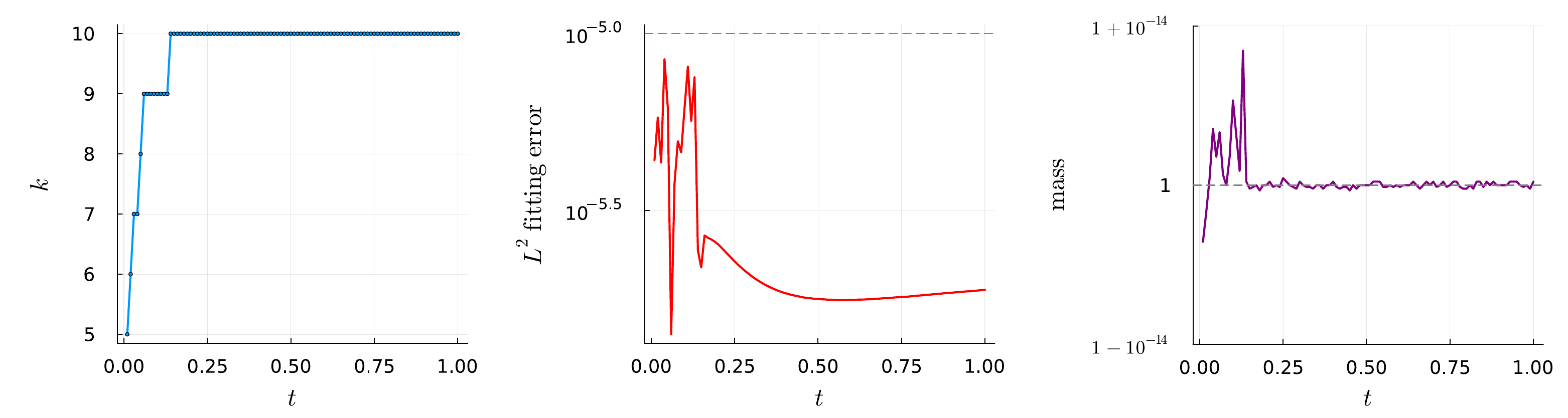}
  \caption{Diagnostics for the double-well problem with $\tau=10^{-3}$. From left to right: number of Gaussian components $k$, transport fitting error, total mass.}\label{fig:dw_diag}
\end{figure}
Figure~\ref{fig:dw_snapshots} shows snapshots of the LCG approximation computed with $\tau=10^{-3}$. The approximation exhibits pronounced symmetry for $t\geq0.2$, although no symmetry is explicitly imposed in the method. Figure~\ref{fig:dw_diag} provides further diagnostics for the same computation. The number of Gaussian components initially increases before stabilizing at $k=10$ around $t=0.2$. The fitting errors in all transport steps remain below the prescribed tolerance of $10^{-5}$, indicating that the Levenberg--Marquardt-based fitting procedure consistently achieves the required accuracy. Mass is conserved to machine precision throughout the simulation. We also report that the total runtime for the $1000$ time steps is $4.59$ seconds on our device.

\begin{figure}[tbp]
  \centering
  \includegraphics[width=0.98\textwidth]{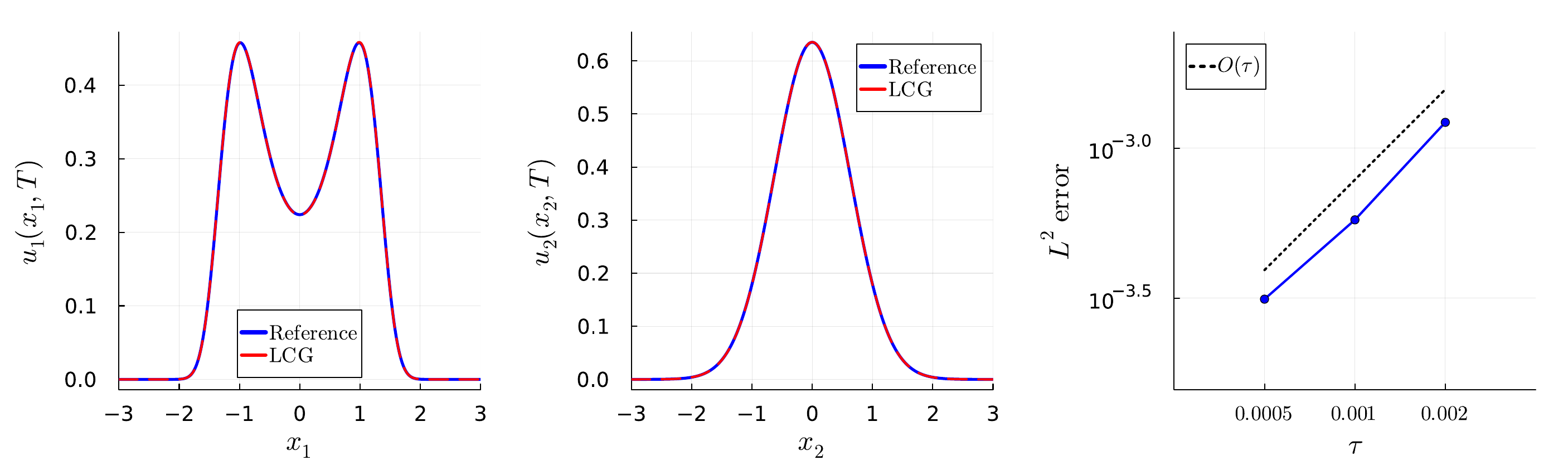}
  \caption{Left and middle: marginal densities for the double-well problem at $T=1$, computed with $\tau=10^{-3}$ (red dashed lines) and compared with the reference solution (blue solid lines). Right: absolute $L^2$ error versus the time step $\tau$ (blue dotted line), with the black dotted line indicating first-order scaling.}\label{fig:dw_marginals}
\end{figure}
The left and middle panels of Figure~\ref{fig:dw_marginals} show the marginal densities at the final time, both of which agree closely with the reference solution. According to Theorem~\ref{thm:global_error}, choosing the fitting tolerance $\varepsilon_{\mathrm{tol}}=\mathcal O(\tau^2)$ yields first-order convergence in time. To examine this behavior numerically, we compute the final-time errors for $\tau=5\times10^{-4}$, $10^{-3}$, and $2\times10^{-3}$. As shown in the right panel of Figure~\ref{fig:dw_marginals}, the errors exhibit approximately first-order decay as $\tau$ decreases, consistent with the theoretical estimate and the order of Lie--Trotter splitting.

\begin{figure}[tbp]
  \centering
  \includegraphics[width=0.55\textwidth]{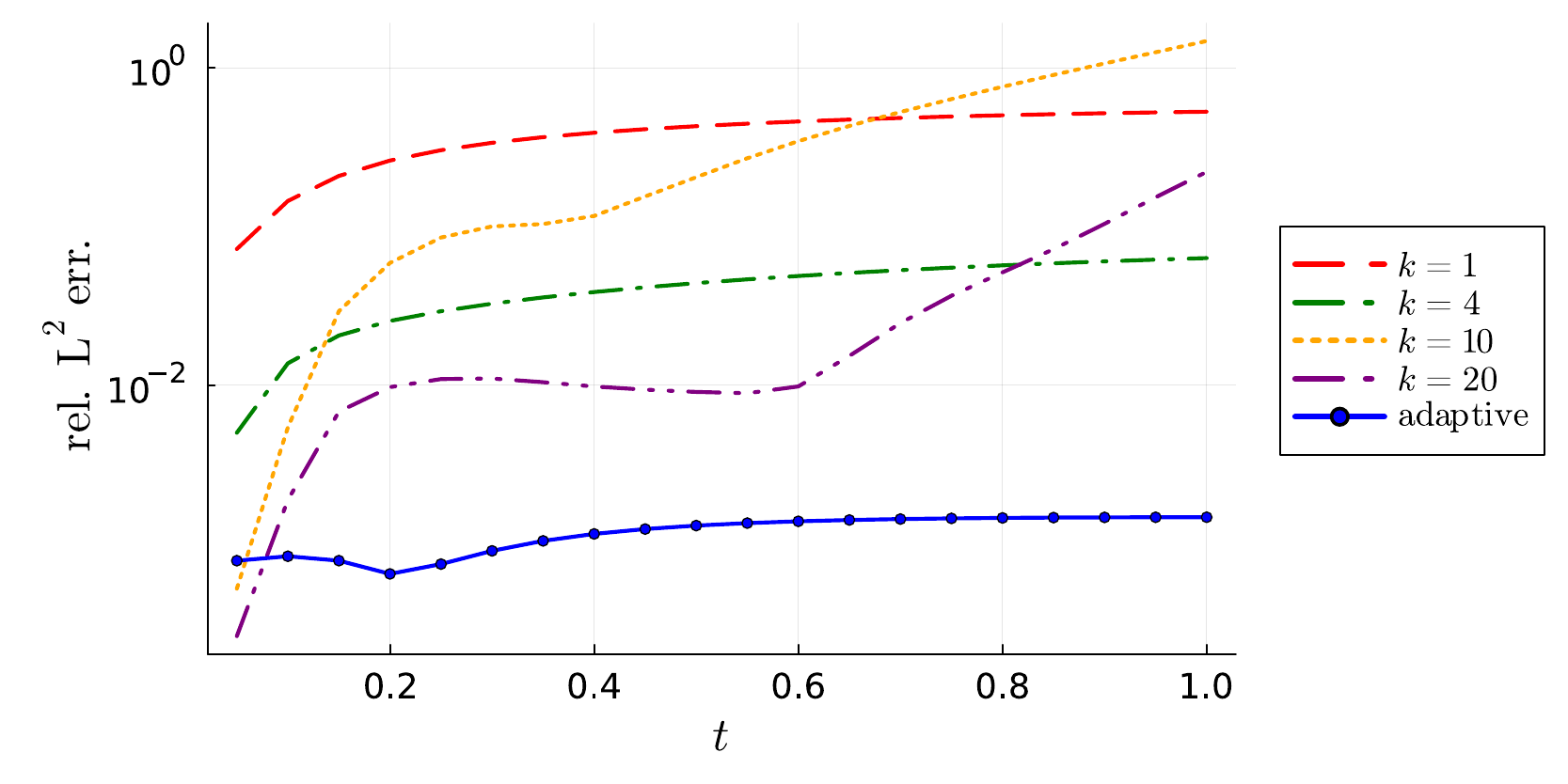}
  \caption{Relative $L^2$ error for the double-well problem with $\tau=10^{-3}$: adaptive scheme and fixed numbers of Gaussian components $k=1,4,10,20$.}\label{fig:dw_fixedK}
\end{figure}
Figure~\ref{fig:dw_fixedK} compares the adaptive scheme with fixed-$k$ runs. The adaptive error remains essentially constant over the whole time interval. With few components, the fixed-$k$ errors are much larger from the start. With as many components as the adaptive run, or twice as many, the fixed-$k$ runs are comparable to the adaptive scheme at early times, but their errors grow rapidly and end well above the adaptive error. Thus, even with the same or twice the maximal number of components, a fixed-$k$ representation does not reproduce the accuracy of the adaptive scheme, and its fitting errors exceed the prescribed tolerance. Unlike in the adaptive scheme, components can neither be introduced where new structure forms nor removed once they become redundant.

\subsection{Two-dimensional radial drift}\label{sec:num_radial}
In this example, we study the formation of an annular density and compare two strategies for adapting the Gaussian representation. We consider the radial drift field
\begin{equation*}
  F(x)=4(1-|x|^2)x, \quad x\in\R^2,
\end{equation*}
with affine part $F_{\degone}=4x$ and nonlinear residual $F_{\deghigher}=-4|x|^2x$. We set the diffusion coefficient to $\kappa=1$ and take the initial density to be the single Gaussian $u_0(x)=\frac{8}{\pi}\exp(-8|x|^2)$. Since the drift is radial and the initial density is radially symmetric, the solution remains radially symmetric. The corresponding Fokker--Planck equation therefore reduces to the following one-dimensional radial equation
\begin{equation}\label{eq:radial_pde}
  \partial_t u=-\frac1r\partial_r\bigl[r(4r-4r^3)u\bigr]
  +\kappa\left(u_{rr}+\frac{u_r}{r}\right),\qquad r>0.
\end{equation}
Here, $r=|x|$ is the radial coordinate, and $u(r,t)$ denotes the density per unit area, satisfying the regularity condition $u_r(0,t)=0$ and the mass normalization $2\pi\int_0^\infty u(r,t)r\,\d r=1$.

We apply the proposed method in two different ways. The Cartesian LCG solves the full two-dimensional equation using unrestricted Gaussian components, whereas the radial LCG solves \eqref{eq:radial_pde} using centered isotropic components of the form $\exp(-\frac{a_i^2}{2}r^2)$. Intuitively, the centered isotropic LCG representation may offer improved scalability to higher dimensions. Both computations are performed up to $T=1$ with $\tau=10^{-3}$, fitting tolerance $\varepsilon_{\mathrm{tol}}=5\times10^{-4}$, and squared recompression tolerance $\delta^2=10^{-7}$. The reference solution is obtained by solving \eqref{eq:radial_pde} on $[0,3]$ using a conservative finite-volume scheme with Chang--Cooper fluxes \cite{CHANG19701} adapted to radial geometry. We use $1000$ uniform radial cells, impose zero-flux boundary conditions, and advance the solution using Crank--Nicolson time stepping with a time step of $2.5\times10^{-5}$.

For the radial LCG ansatz, the refinement strategy described in Section~\ref{sec:component_splitting} is adapted to preserve radial symmetry. Radial moments are evaluated as centered isotropic two-dimensional Gaussian integrals with respect to the radial measure $2\pi r\,\d r$. Mass-constrained fitting and recompression otherwise follow the same procedures as in the Cartesian computation.

\begin{figure}[tbp]
  \centering
  \includegraphics[width=0.98\textwidth]{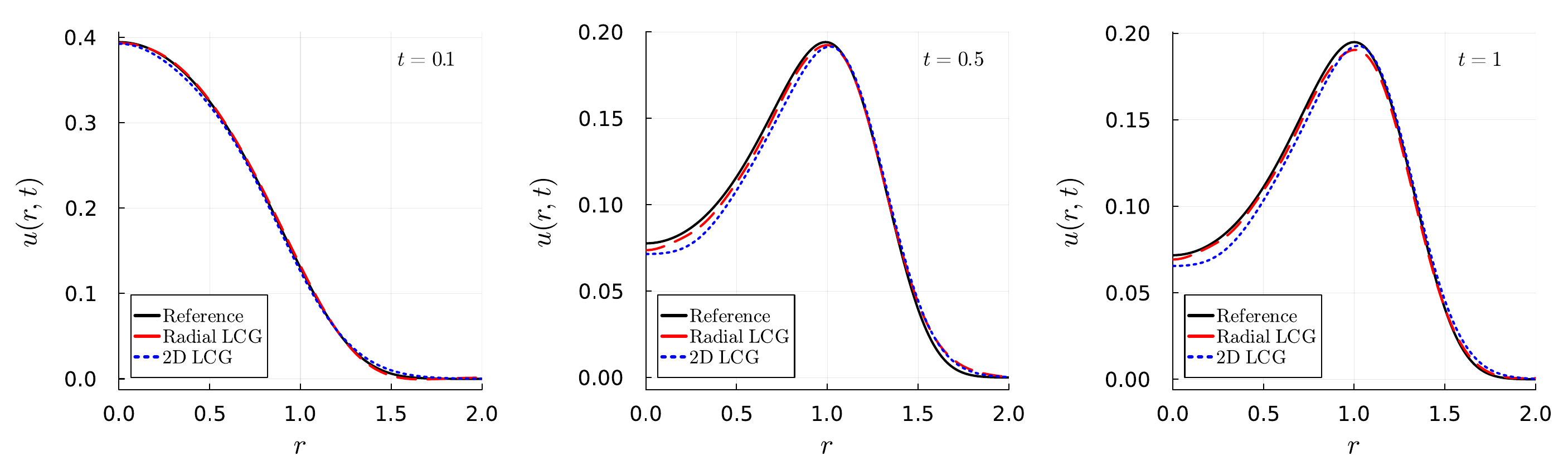}
  \caption{Radial profiles computed with $\tau=10^{-3}$. From left to right: $t=0.1$, $0.5$, and $1$. The radial LCG (red dashed lines) and angular average of the Cartesian LCG (blue dotted lines) are compared with the reference solution (black solid lines).}\label{fig:rad_profiles}
\end{figure}
\begin{figure}[tbp]
  \centering
  \includegraphics[width=0.98\textwidth]{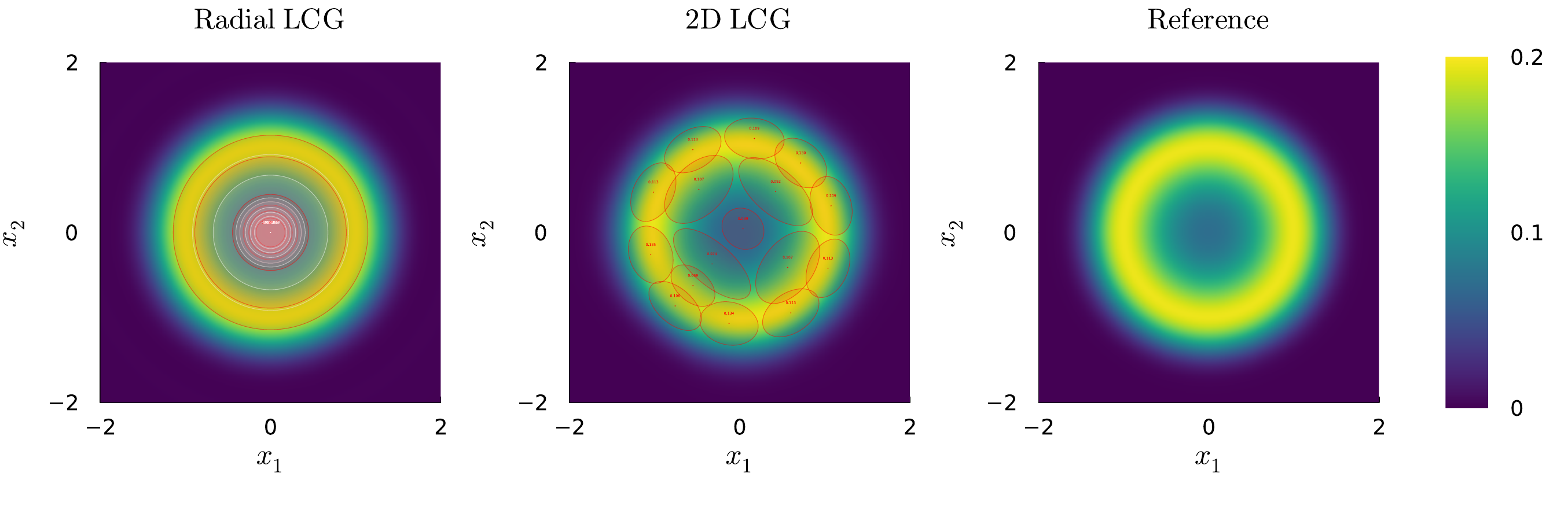}
  \caption{Final densities for the radial problem at $T=1$, computed with $\tau=10^{-3}$. From left to right: radial LCG; Cartesian LCG; reference solution. Red and white ellipses represent Gaussian components with positive and negative coefficients, respectively.}\label{fig:rad_final}
\end{figure}
Figure~\ref{fig:rad_profiles} shows the transition from a centrally peaked density to an annular profile. Both approximations capture this transition and agree closely with the reference solution, with the radial LCG providing slightly better agreement near $r=0$. The radial profiles of the Cartesian LCG are obtained by angular averaging over $128$ equally spaced angles. Figure~\ref{fig:rad_final} compares the final densities in the plane. The radial LCG preserves radial symmetry by construction, whereas the Cartesian LCG forms the ring using displaced components and exhibits a small angular variation.

\begin{table}[tbp]
  \centering
  \begin{tabular}{l c c c c}
    \hline
    Method        & $k(T)$ & profile error         & time (s) & $\max\sqrt{\mathcal E}/\varepsilon_{\mathrm{tol}}$ \\
    \hline
    Cartesian LCG & 16     & $4.2243\times10^{-2}$ & 17.62    & 0.998                                              \\
    Radial LCG    & 14     & $2.3027\times10^{-2}$ & 20.52    & 2.061                                              \\
    \hline
  \end{tabular}
  \caption{Comparison of the Cartesian and radial LCGs at $T=1$ with $\tau=10^{-3}$.}\label{tab:radial_comparison}
\end{table}
Table~\ref{tab:radial_comparison} compares the radial and Cartesian LCG approximations in terms of the final number of Gaussian components, radial-profile error, runtime, and maximum fitting error. The radial LCG achieves a smaller profile error with fewer components, although it incurs a slightly higher computational cost in these experiments. For the Cartesian LCG, the relative $L^2$ error of the full two-dimensional density, evaluated using polar quadrature, is $4.7709\times 10^{-2}$, which is larger than the corresponding radial-profile error. This difference arises because angular averaging suppresses the deviations from radial symmetry observed in Figure~\ref{fig:rad_final}, which contribute to the error of the full two-dimensional approximation.

\begin{figure}[tbp]
  \centering
  \includegraphics[width=0.98\textwidth]{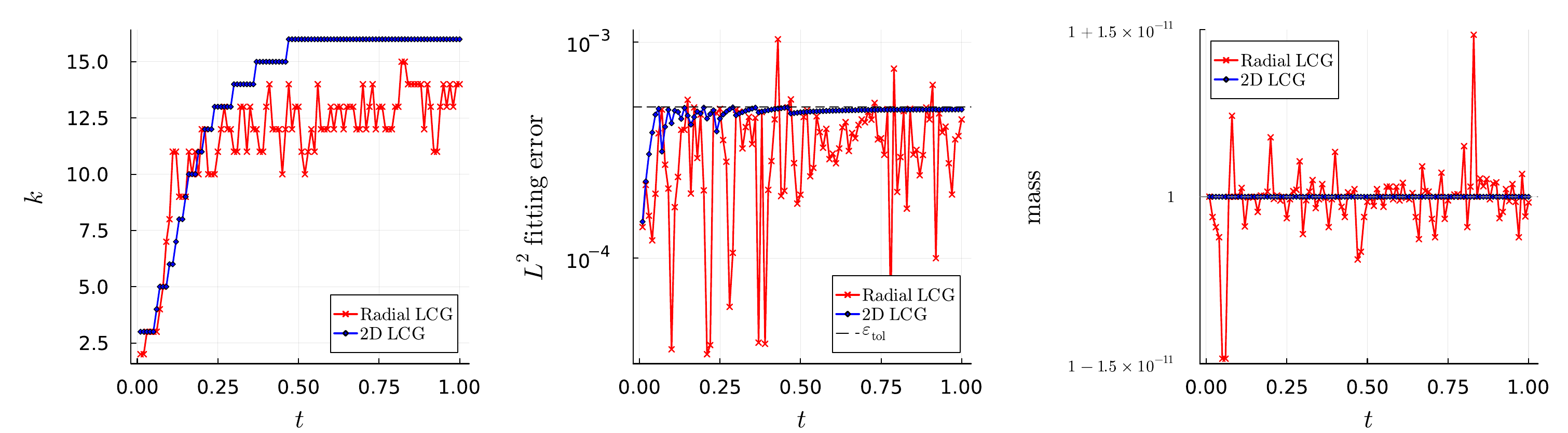}
  \caption{Diagnostics for the radial problem with $\tau=10^{-3}$. From left to right: number of Gaussian components $k$, transport fitting error, total mass.}\label{fig:rad_diag}
\end{figure}
Figure~\ref{fig:rad_diag} presents further diagnostics for both computations. For the Cartesian LCG, the number of components increases to $k=16$ and subsequently remains constant, whereas the radial LCG undergoes further enrichment and recompression, ending with $k=14$ components. All recorded fitting errors for the Cartesian LCG remain below the prescribed tolerance $\varepsilon_{\mathrm{tol}}$. In contrast, the radial LCG exceeds this tolerance at six time steps, with a maximum fitting error of $2.061\varepsilon_{\mathrm{tol}}$. Despite these occasional tolerance violations, both approximations exhibit excellent mass conservation, with mass deviations bounded by $1.5\times10^{-11}$ for the radial LCG and $6\times10^{-16}$ for the Cartesian LCG.

\begin{figure}[tbp]
  \centering
  \includegraphics[width=0.6\textwidth]{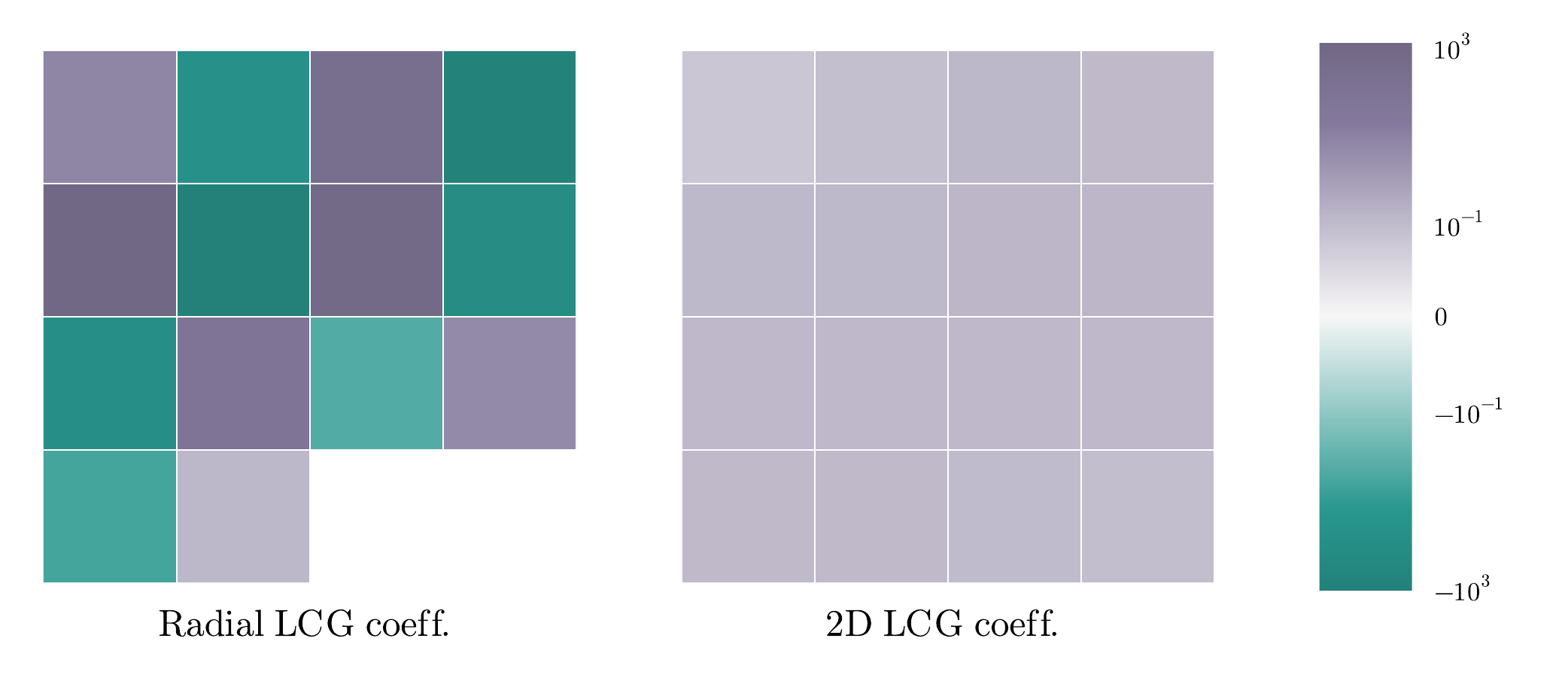}
  \caption{Gaussian coefficients for the radial problem at $T=1$, computed with $\tau=10^{-3}$. Left: radial LCG. Right: Cartesian LCG. Each cell represents one coefficient, and coefficients are ordered by increasing $\operatorname{tr}(\Sigma_i)$.}\label{fig:rad_coefficients}
\end{figure}
Finally, we discuss the numerical stability of the radial LCG representation. We visualize the Gaussian coefficients of the final solution in Figure~\ref{fig:rad_coefficients}. The radial coefficients range from approximately $-1.06\times10^3$ to $1.15\times10^3$, whereas the Cartesian coefficients lie between $0.039$ and $0.135$. Moreover, all Cartesian coefficients are positive, while seven of the fourteen radial coefficients are negative. The radial representation thus relies on substantial cancellation between components, making fitting and residual evaluation sensitive to roundoff errors. The need for signed coefficients to reproduce an annular profile within the centered isotropic ansatz can be understood by expressing the Gaussian approximation as a discrete Laplace transform in the squared radial variable. Consider the approximation
$$
  f(r)\approx f_{\mathrm{LCG}}(r)=\sum_i c_i\exp\left(-\frac{a_i^2}{2}r^2\right).
$$
Setting $r=\sqrt{s}$ and $q_i=a_i^2/2$, we obtain the discrete Laplace representation
$$
  g_{\mathrm{LCG}}(s)=f_{\mathrm{LCG}}(\sqrt{s})
  =\sum_i c_i e^{-sq_i}
  =\int_0^\infty e^{-st}\,\mathrm d\mu(t),
  \quad\mu=\sum_i c_i\delta_{q_i}.
$$
Here, $\delta_{q_i}$ denotes the Dirac measure concentrated at $q_i$. If all coefficients are nonnegative, then
$$
  g_{\mathrm{LCG}}'(s)=-\sum_i c_i q_i e^{-sq_i}\leq0,
$$
so both $g_{\mathrm{LCG}}$ and $f_{\mathrm{LCG}}$ are nonincreasing. Such a representation cannot reproduce an annular density whose radial profile increases away from the origin before attaining a maximum at a positive radius. Signed coefficients are therefore necessary to capture this structure using centered isotropic Gaussian components.

\subsection{Two-dimensional nonlinear oscillator}\label{sec:num_duffing}

This example tests the method on a nonlinear oscillator whose drift couples the two coordinates and deforms the density through stretching, shear, and rotation. We consider the Duffing-type drift field
\begin{equation*}
  F(x_1,x_2)=\begin{pmatrix}2x_2\\2x_1-0.8x_2-0.2x_1^3\end{pmatrix},
\end{equation*}
with affine part $F_{\degone}=(2x_2,2x_1-0.8x_2)^\top$ and nonlinear residual $F_{\deghigher}=(0,-0.2x_1^3)^\top$. We set the diffusion coefficient to $\kappa=0.3$ and take the initial density to be the single Gaussian $u_0(x)=\frac{1}{2\pi}\exp(-\frac12|x|^2)$. The solution is computed up to $T=1$ with $\tau=10^{-3}$, the fitting tolerance $\varepsilon_{\mathrm{tol}}=5\times10^{-5}$, and the squared recompression tolerance $\delta^2=\frac12\varepsilon_{\mathrm{tol}}^2$. The reference solution is computed on $[-6,6]^2$ using Chang--Cooper finite-volume fluxes \cite{CHANG19701} on a uniform $401\times401$ grid, zero-flux boundary conditions, and symmetric directional Crank--Nicolson sweeps with time step size $2.5\times10^{-4}$.

\begin{figure}[tbp]
  \centering
  \includegraphics[width=0.98\textwidth]{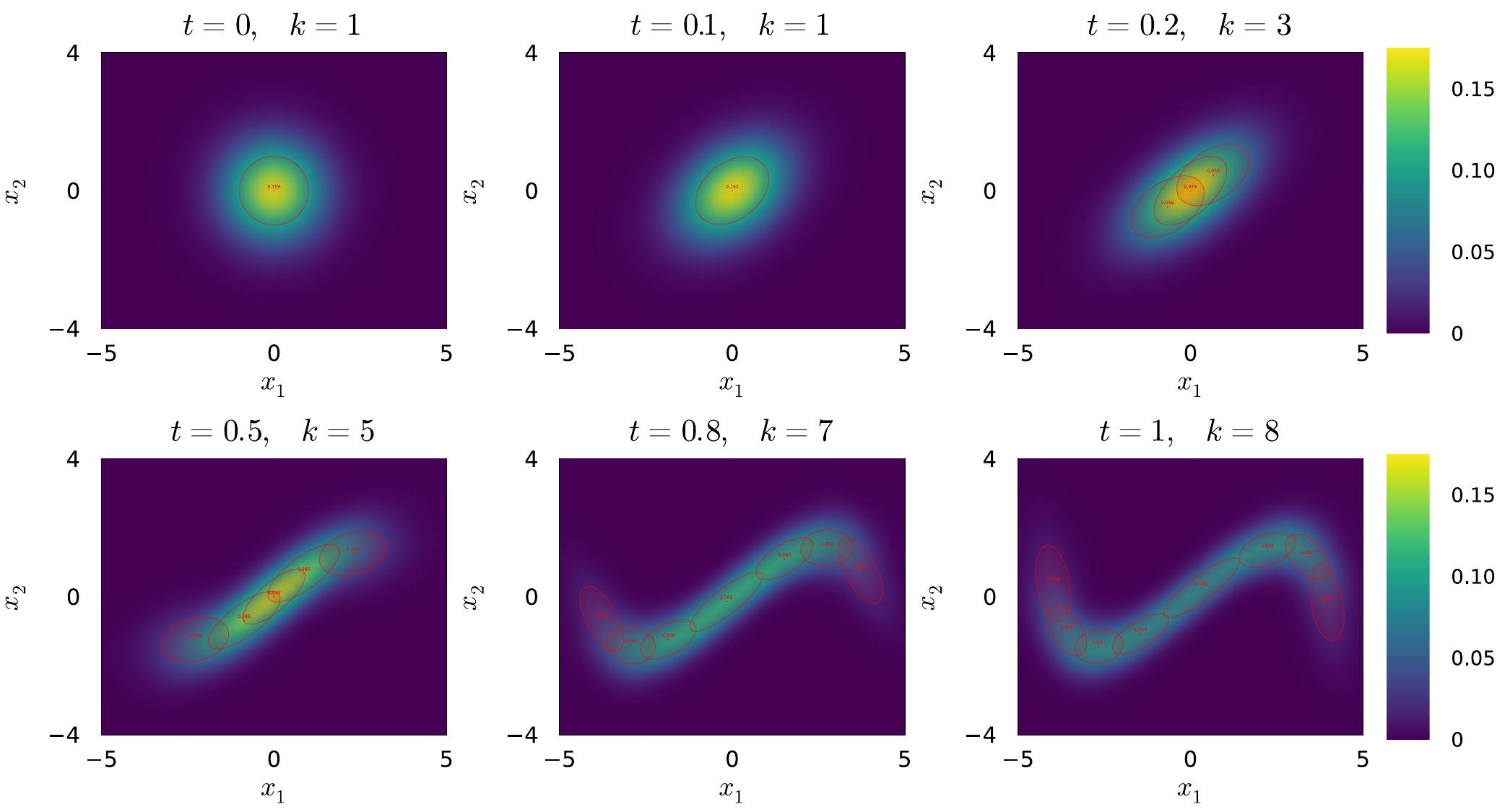}
  \caption{Snapshots of the LCG solution for the nonlinear oscillator with $\tau=10^{-3}$. Red ellipses represent Gaussian components with positive coefficients.}\label{fig:duff_snapshots}
\end{figure}
\begin{figure}[tbp]
  \centering
  \includegraphics[width=0.98\textwidth]{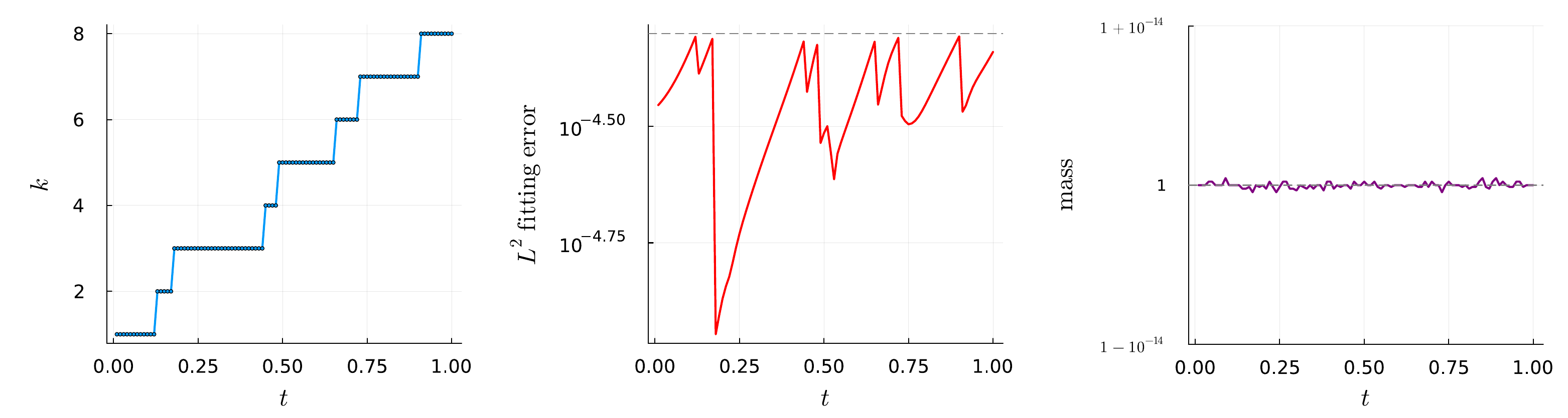}
  \caption{Diagnostics for the nonlinear oscillator with $\tau=10^{-3}$. From left to right: number of Gaussian components $k$, transport fitting error, total mass.}\label{fig:duff_diag}
\end{figure}
Figure~\ref{fig:duff_snapshots} shows snapshots of the LCG approximation computed with $\tau=10^{-3}$. The initially circular Gaussian stretches and develops a curved profile, which is represented by progressively adding and reorienting Gaussian components. Figure~\ref{fig:duff_diag} provides further diagnostics for the same computation. The number of Gaussian components increases from $k=1$ to $k=8$ by $T=1$. The recorded fitting errors remain below the prescribed tolerance of $5\times10^{-5}$, with a maximum of $0.987\varepsilon_{\mathrm{tol}}$. Mass is conserved to machine precision throughout the simulation. The total runtime for the $1000$ time steps is $1.22$ seconds on our device.

\begin{figure}[tbp]
  \centering
  \includegraphics[width=0.98\textwidth]{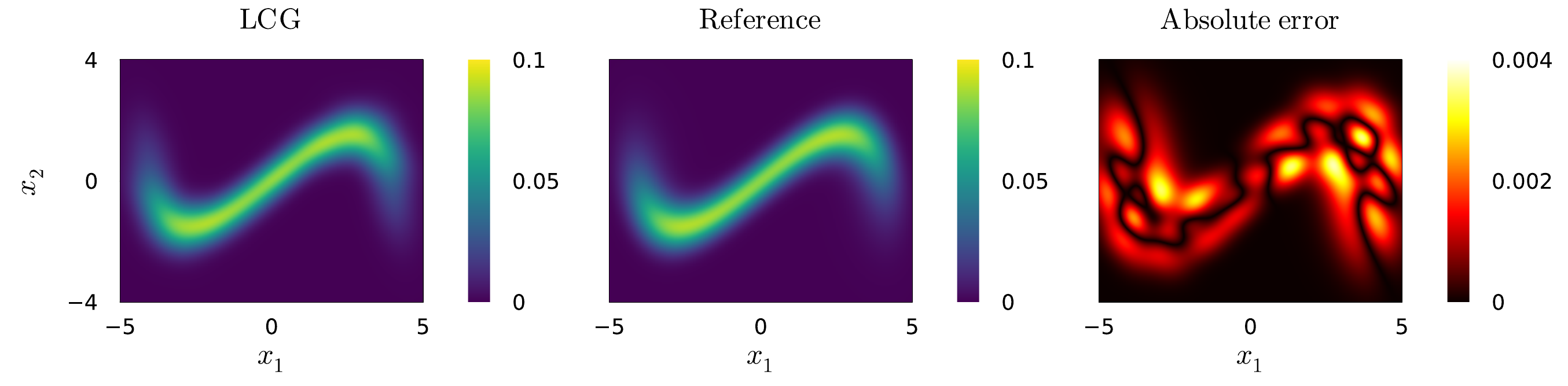}
  \caption{Final densities for the nonlinear oscillator at $T=1$, computed with $\tau=10^{-3}$. From left to right: LCG solution, reference solution, pointwise absolute error.}\label{fig:duff_adi}
\end{figure}
The left and middle panels of Figure~\ref{fig:duff_adi} show the final LCG approximation and the reference solution, respectively, which are in close agreement. The LCG approximation captures the curved structure of the density, with an absolute $L^2$ error of $8.4828\times10^{-3}$ and a relative $L^2$ error of $3.7083\times10^{-2}$. The right panel displays the corresponding pointwise absolute error.

\begin{figure}[tbp]
  \centering
  \includegraphics[width=0.55\textwidth]{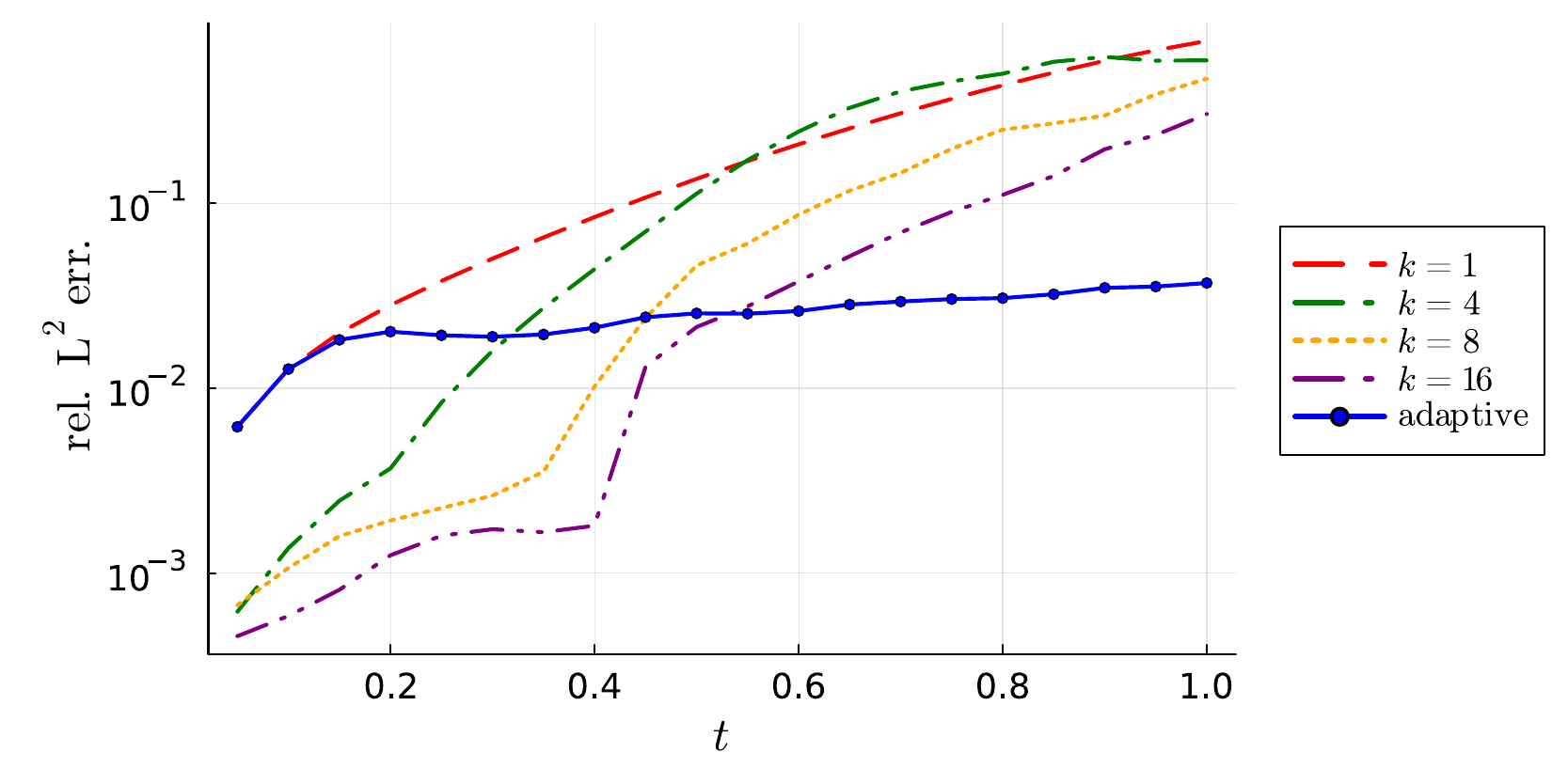}
  \caption{Relative $L^2$ error for the nonlinear oscillator with $\tau=10^{-3}$: adaptive scheme and fixed numbers of Gaussian components $k=1,4,8,16$.}\label{fig:duff_fixedK}
\end{figure}
Figure~\ref{fig:duff_fixedK} shows the same comparison for the nonlinear oscillator. Initially, a single Gaussian satisfies the fitting tolerance, so the adaptive run does not split and coincides with the run with $k=1$. In this phase, the fixed-$k$ runs with more components are more accurate. Once the density develops its curved structure, however, their errors grow steadily, while the adaptive scheme adds components and keeps the error at a nearly constant level. At the final time, all fixed-$k$ runs are clearly less accurate than the adaptive scheme.

\subsection{Three-dimensional damped spinning top}\label{sec:num_spinning_top}
The last example combines nonlinear rotation with linear damping in three dimensions. It tests whether enrichment captures the initial deformation and recompression reduces the representation as the density contracts. We consider a symmetric damped spinning top, whose drift field is
\begin{equation*}
  F(x)=\begin{pmatrix}
    -0.2x_1-x_2x_3 \\
    -0.2x_2+x_3x_1 \\
    -0.1x_3
  \end{pmatrix}.
\end{equation*}
The affine part is $F_{\degone}=(-0.2x_1,-0.2x_2,-0.1x_3)^\top$, and the nonlinear residual contains the two quadratic rotation terms. We set the diffusion coefficient to $\kappa=10^{-4}$ and take the initial density to be the single Gaussian
\begin{equation*}
  u_0(x)=\frac{1}{(2\pi)^{3/2}\,0.05^2\cdot0.4}
  \exp\left(-\frac{(x_1-0.1)^2+(x_2-0.1)^2}{2\cdot0.05^2}
  -\frac{(x_3-1)^2}{2\cdot0.4^2}\right).
\end{equation*}
We compute the solution up to $T=100$ using the time step $\tau=10^{-2}$, the fitting tolerance $\varepsilon_{\mathrm{tol}}=10^{-3}$, and the squared recompression tolerance $\delta^2=\varepsilon_{\mathrm{tol}}^2/2$. In this computation, the fitting objective is normalized by the squared $L^2$ norm of the transport target. To construct the reference solution, we exploit the conditional Gaussian structure of the symmetric top: the third coordinate follows an Ornstein--Uhlenbeck process and is jointly Gaussian with the accumulated rotation angle, allowing density evaluation to be reduced, by conditioning, to a one-dimensional Gauss--Hermite quadrature. We use $96$ quadrature nodes and evaluate the reference density on adaptive $81\times81\times101$ Cartesian grids at the saved times. Further increases in the number of quadrature nodes produce negligible changes in the reference.

\begin{figure}[tbp]
  \centering
  \includegraphics[width=0.98\textwidth]{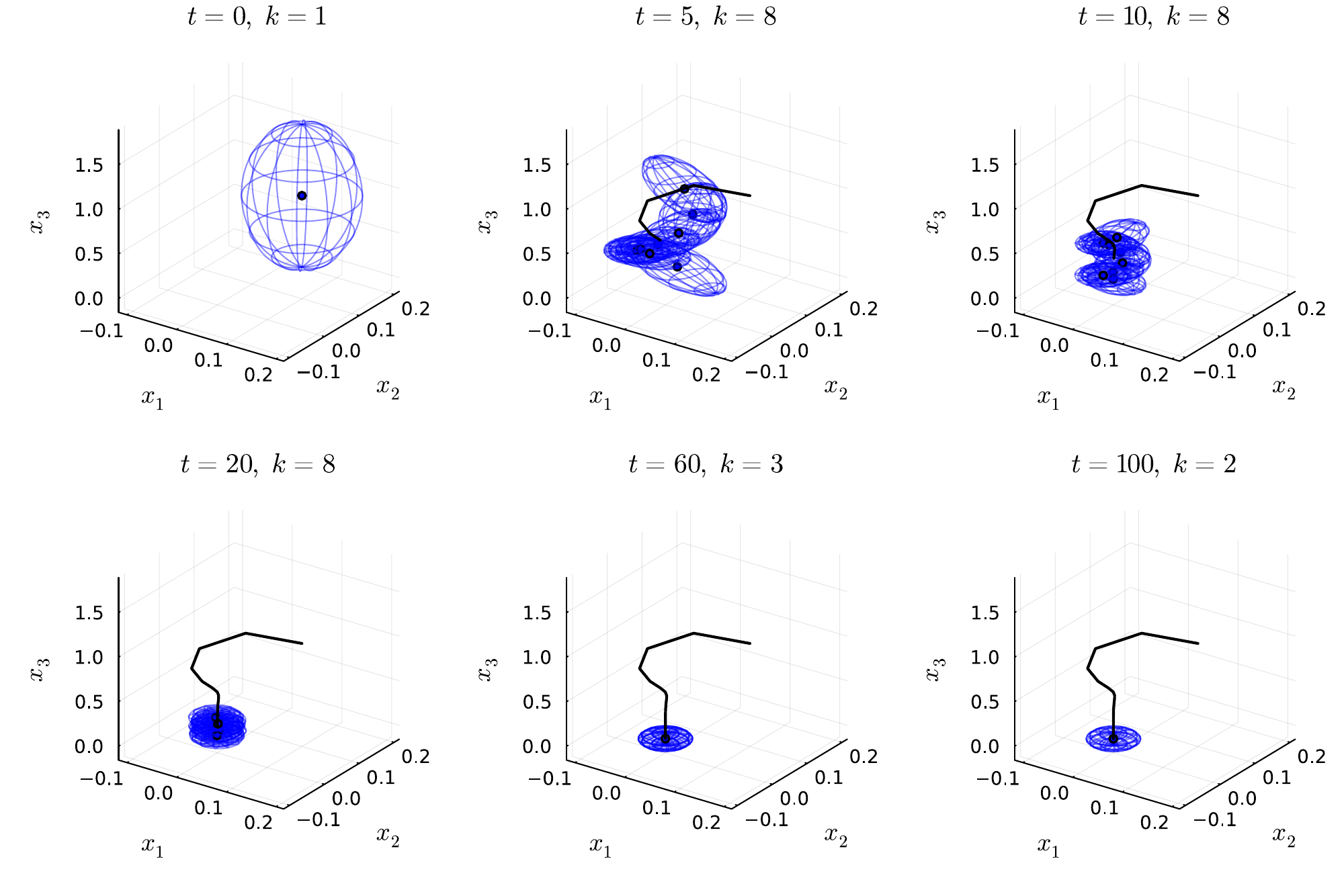}
  \caption{Snapshots of the Gaussian components for the damped spinning top with $\tau=10^{-2}$. Blue ellipsoids represent components with positive coefficients. The black line shows the trajectory of the LCG mean.}\label{fig:top_traj}
\end{figure}
\begin{figure}[htb!]
  \centering
  \includegraphics[width=0.98\textwidth]{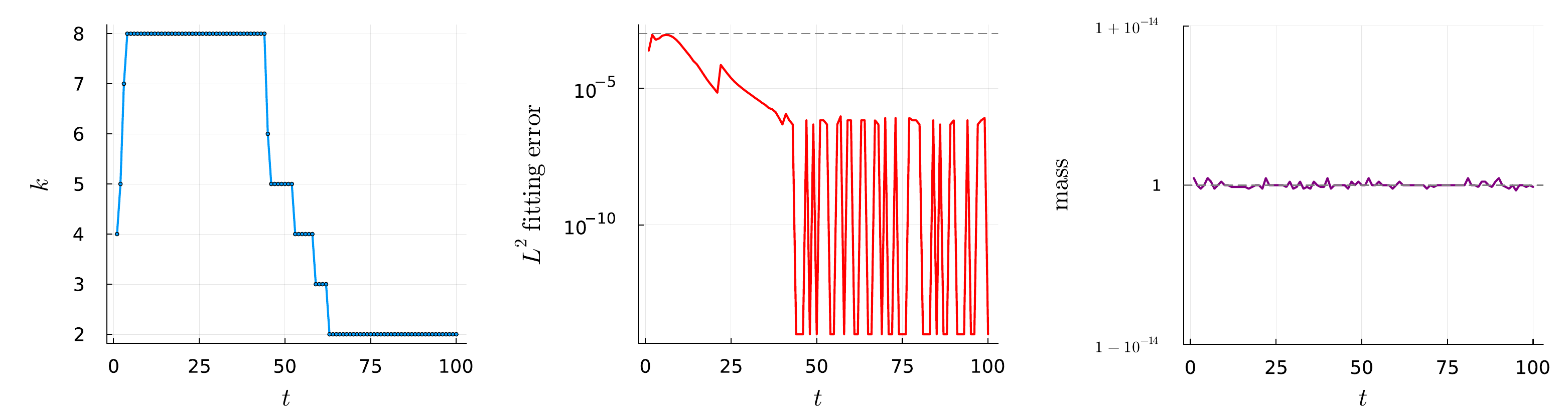}
  \caption{Diagnostics for the damped spinning top with $\tau=10^{-2}$. From left to right: number of Gaussian components $k$, transport fitting error, total mass.}\label{fig:top_diag}
\end{figure}
Figure~\ref{fig:top_traj} shows snapshots of the Gaussian components computed with $\tau=10^{-2}$, together with the trajectory of the LCG mean. The components initially spread along the rotational motion before contracting toward the origin under damping. Further diagnostics are presented in Figure~\ref{fig:top_diag}. The number of Gaussian components increases from $k=1$ to $k=8$ during the initial transient and subsequently decreases to $k=2$ by $T=100$. All recorded fitting errors remain below the prescribed tolerance of $10^{-3}$, with a maximum of $0.919\varepsilon_{\mathrm{tol}}$, while mass is conserved to machine precision throughout the simulation. The total runtime for the $10000$ time steps is $40.34$ seconds on our device.

\begin{figure}[tbp]
  \centering
  \includegraphics[width=0.8\textwidth]{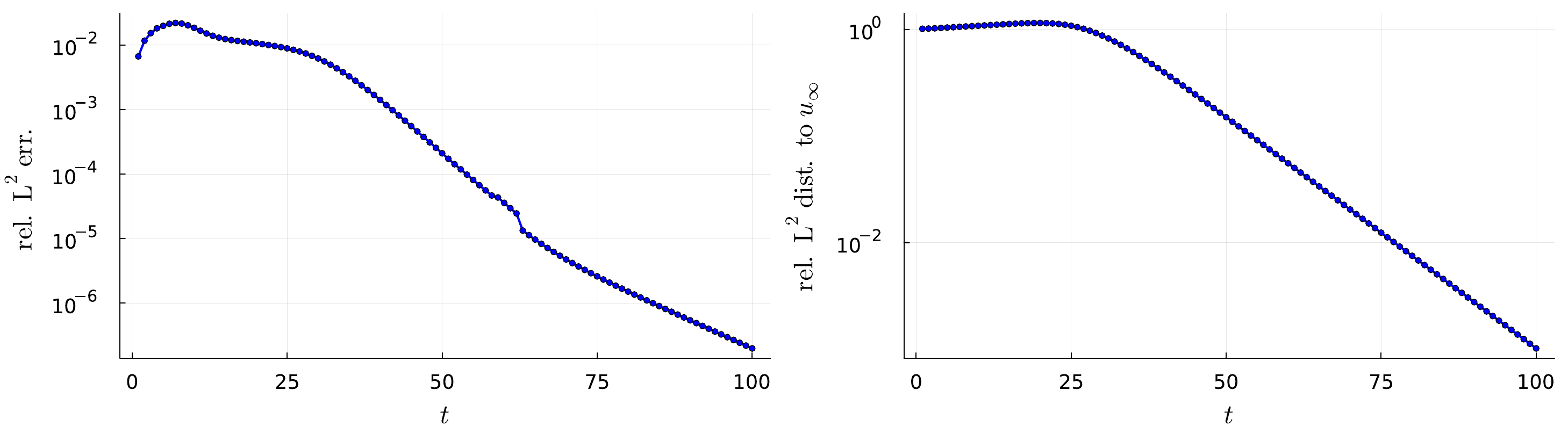}
  \caption{Error evolution for the damped spinning top, computed with $\tau=10^{-2}$. Left: relative $L^2$ error with respect to the reference solution. Right: relative $L^2$ distance $\norm{u^n-u_\infty}_{L^2}/\norm{u_\infty}_{L^2}$ to the equilibrium $u_\infty$.}\label{fig:top_convergence_pair}
\end{figure}
Figure~\ref{fig:top_convergence_pair} shows the relative $L^2$ error with respect to the reference solution (left) and the relative $L^2$ distance to the equilibrium (right),
\begin{equation*}
  u_\infty(x)=\frac{1}{(2\pi)^{3/2}\sqrt{250\kappa^3}}
  \exp\left(-\frac{x_1^2+x_2^2}{10\kappa}-\frac{x_3^2}{20\kappa}\right),
\end{equation*}
which is the stationary solution for this drift. The LCG approximation remains close to the reference solution and converges to the equilibrium.

\subsection{Discussion}\label{sec:num_discussion}
These examples demonstrate that the proposed method captures the formation of multiple peaks, annular profiles, and curved densities through an adaptive Gaussian representation without a spatial grid. Enrichment introduces additional components as the density deforms, while recompression removes redundant components as its structure simplifies. The comparison with fixed-$k$ runs shows that this adaptivity is essential: with the same or twice the maximal number of components of the adaptive run, the fixed-$k$ runs end with considerably larger errors. Mass deviations remain below $1.5\times10^{-11}$ in all four examples. The radial comparison also highlights the benefits and numerical challenges of imposing symmetry: centered isotropic Gaussians yield a smaller profile error with fewer components, but the resulting representation requires signed coefficients, whose cancellation makes accurate residual evaluation more demanding.

\section{Conclusion}\label{sec:conclusion}
We develop an adaptive, mesh-free method for linear Fokker--Planck equations with polynomial drift on $\R^d$. Based on Lie--Trotter splitting, the method combines exact evolution under the affine Fokker--Planck operator with mass-constrained $L^2$ fitting of a Taylor approximation to the density transported by the nonlinear drift. Gaussian moment formulas yield an analytically computable fitting error, which determines whether enrichment is required, while recompression controls redundancy in the Gaussian representation. Under suitable stability and regularity assumptions, the local and global error estimates quantify the contributions from splitting, Taylor truncation, fitting, and recompression. Numerical experiments in two and three dimensions illustrate the method's ability to resolve multimodal, annular, and anisotropic densities while adapting the number of Gaussian components throughout the evolution.

An important direction for future work is to improve scalability to substantially higher dimensions, with particular attention to the computational cost of Gaussian fitting and moment evaluation. Extensions to nonlinear Fokker--Planck equations and higher-order splitting schemes are also of interest.

\appendix
\section{Proofs and Taylor coefficients for the transport projection}\label{sec:transport_taylor_appendix}

\subsection{Proof of Lemma~\ref{lem:transport_cross}}
\begin{proof}
  Changing variables $y=\psi_t(x)$ in the conservative representation gives
  \begin{equation}\label{eq:Icross_cov}
    I_{\mathrm{\times}}(t)=\int_{\R^d}u_{c,w}(\psi_t(x))v_0(x)\,\d x.
  \end{equation}
  Set $L=F_{\deghigher }\cdot\nabla$. Repeated differentiation yields
  $\partial_t^\ell g_{w_i}(\psi_t(x))=(L^\ell g_{w_i})(\psi_t(x))$.
  Define
  $$
    p_{\times,\ell,i}(x)=\frac{L^\ell g_{w_i}(x)}{\ell!\,g_{w_i}(x)},\qquad
    p_{\times,i}^{(r)}(x;t)=\sum_{\ell=0}^r t^\ell p_{\times,\ell,i}(x).
  $$
  Writing $P_i=U_i^\top U_i$, these coefficients satisfy
  $$
    p_{\times,0,i}=1,\qquad
    p_{\times,\ell+1,i}=\frac{F_{\deghigher }\cdot\nabla p_{\times,\ell,i}
      -\bigl(F_{\deghigher }^\top P_i(x-m_i)\bigr)p_{\times,\ell,i}}{\ell+1}.
  $$
  In particular, all these coefficients are polynomial. Since $\mathcal B$ is the formal adjoint of $L$, repeated integration by parts in \eqref{eq:b_vector_transport} gives
  $$
    \begin{aligned}
      b_i(w;t)=\sum_{\ell=0}^r\frac{t^\ell}{\ell!}
      \langle L^\ell g_{w_i},v_0\rangle=\sum_{j=1}^{k_0}c_{0,j}
      \int_{\R^d}p_{\times,i}^{(r)}(x;t)g_{w_i}(x)g_{w_{0,j}}(x)\,\d x.
    \end{aligned}
  $$
  The boundary terms vanish because all functions involved are polynomials times Gaussians. Taylor's formula gives, for each fixed $x$,
  \begin{equation}\label{eq:generic_cross_taylor}
    g_{w_i}(\psi_t(x))=g_{w_i}(x)
    \left(p_{\times,i}^{(r)}(x;t)+\mathcal O(t^{r+1})\right).
  \end{equation}
  For the integrated estimate, use the absolute remainder bound
  $$
    \left|g_{w_i}(\psi_t(x))-g_{w_i}(x)p_{\times,i}^{(r)}(x;t)\right|
    \leq\frac{t^{r+1}}{(r+1)!}\|L^{r+1}g_{w_i}\|_{L^\infty}.
  $$
  The norm on the right is finite because $L^{r+1}g_{w_i}$ is a polynomial times a Gaussian. Substitution into \eqref{eq:Icross_cov} proves the claim with
  $$
    C_{\times}=\frac{\|v_0\|_{L^1}}{(r+1)!}
    \sum_{i=1}^k|c_i|\|L^{r+1}g_{w_i}\|_{L^\infty}.
  $$

  For completeness, we give explicit cross coefficients for the choice $r=3$ used in Section~\ref{sec:numerical_examples}. Here $DF_{\deghigher }$ is the Jacobian matrix. The vector $\D^2F_{\deghigher }(v,w)$ has components $v^\top \D^2(F_{\deghigher })_\alpha w$.
  Let
  \[
    A(x)=\D F_{\deghigher }(x)F_{\deghigher }(x), \qquad
    B(x)=\D^2F_{\deghigher }(x)(F_{\deghigher }(x),F_{\deghigher }(x))+\D F_{\deghigher }(x)A(x).
  \]
  Then the flow generated by \eqref{eq:charac} satisfies
  \begin{equation}\label{eq:app_flow_expansion}
    \psi_t(x)=x+tF_{\deghigher }(x)+\frac{t^2}{2}A(x)+\frac{t^3}{6}B(x)+\mathcal O(t^4).
  \end{equation}
  For $w_i=(m_i,U_i)$, set $P_i=U_i^\top U_i$. Substituting \eqref{eq:app_flow_expansion} into $g_{w_i}(\psi_t(x))$ gives
  $$
    g_{w_i}(\psi_t(x))=g_{w_i}(x)\exp\left(
    t\theta_{1,i}(x)+t^2\theta_{2,i}(x)+t^3\theta_{3,i}(x)+\mathcal O(t^4)
    \right),
  $$
  where
  $$
    \theta_{1,i}(x)=-F_{\deghigher }(x)^\top P_i(x-m_i),
  $$
  $$
    \theta_{2,i}(x)=-\frac12 A(x)^\top P_i(x-m_i)-\frac12 F_{\deghigher }(x)^\top P_iF_{\deghigher }(x),
  $$
  and
  $$
    \theta_{3,i}(x)=-\frac16 B(x)^\top P_i(x-m_i)-\frac12 F_{\deghigher }(x)^\top P_iA(x).
  $$
  The third-order truncated cross polynomial is therefore
  \begin{equation*}
    p_{\times,i}^{(3)}(x;t)=1+t\theta_{1,i}(x)+t^2\left(\theta_{2,i}(x)+\frac12\theta_{1,i}(x)^2\right)
    +t^3\left(\theta_{3,i}(x)+\theta_{1,i}(x)\theta_{2,i}(x)+\frac16\theta_{1,i}(x)^3\right),
  \end{equation*}
  which is just the Taylor expansion of the exponential up to order three.
  Thus
  $$
    g_{w_i}(\psi_t(x))=g_{w_i}(x)\left(p_{\times,i}^{(3)}(x;t)+\mathcal O(t^4)\right).
  $$
  In the notation of \eqref{eq:generic_cross_taylor},
  $$
    p_{\times,0,i}=1,\qquad
    p_{\times,1,i}=\theta_{1,i},\qquad
    p_{\times,2,i}=\theta_{2,i}+\frac12\theta_{1,i}^2,\qquad
    p_{\times,3,i}=\theta_{3,i}+\theta_{1,i}\theta_{2,i}+\frac16\theta_{1,i}^3. \qedhere
  $$
\end{proof}

\subsection{Proof of Lemma~\ref{lem:transport_norm}}
\begin{proof}
  Expanding $\norm{\widetilde u_{\adv}^{(r)}(t)}^2$ using \eqref{eq:transport_taylor_target} gives \eqref{eq:hat_I_nl}. For its evaluation by Gaussian moments, define
  $$
    p_{\adv,\ell,j}(x)=\frac{\mathcal B^\ell g_{w_{0,j}}(x)}{\ell!\,g_{w_{0,j}}(x)},\qquad
    p_{\adv,j}^{(r)}(x;t)=\sum_{\ell=0}^r t^\ell p_{\adv,\ell,j}(x).
  $$
  Write $P_{0,j}=U_{0,j}^\top U_{0,j}$. Applying $\mathcal B$ to a polynomial times the $j$-th input Gaussian gives
  $$
    \begin{aligned}
      \mathcal B(pg_{w_{0,j}})
      =g_{w_{0,j}}\Bigl(\bigl(F_{\deghigher}^\top P_{0,j}(x-m_{0,j})
      -\nabla\cdot F_{\deghigher}\bigr)p-F_{\deghigher}\cdot\nabla p\Bigr).
    \end{aligned}
  $$
  The density coefficients are determined recursively by
  \begin{equation*}
    \begin{aligned}
      p_{\adv,0,j} & =1,                                           \\
      (\ell+1)p_{\adv,\ell+1,j}
                   & =\bigl(F_{\deghigher}^\top P_{0,j}(x-m_{0,j})
      -\nabla\cdot F_{\deghigher}\bigr)p_{\adv,\ell,j}-F_{\deghigher}\cdot\nabla p_{\adv,\ell,j}.
    \end{aligned}
  \end{equation*}
  In particular, all coefficients are polynomial, and substitution into \eqref{eq:transport_taylor_target} gives
  $$
    \widetilde u_{\adv}^{(r)}(t,x)
    =\sum_{j=1}^{k_0}c_{0,j}p_{\adv,j}^{(r)}(x;t)g_{w_{0,j}}(x).
  $$
  Squaring this representation and integrating yields
  $$
    \widehat I_{\adv}(t)=\sum_{i,j=1}^{k_0}c_{0,i}c_{0,j}
    \int_{\R^d}p_{\adv,i}^{(r)}(x;t)p_{\adv,j}^{(r)}(x;t)
    g_{w_{0,i}}(x)g_{w_{0,j}}(x)\,\d x.
  $$
  For an order-$r$ approximation, the recurrence is applied for $\ell=0,\ldots,r-1$, and the full squared norm retains all powers through $t^{2r}$.

  For the remainder estimate, Taylor's formula in $L^2$ gives
  $$
    u_{\adv}(t)-\widetilde u_{\adv}^{(r)}(t)
    =\frac1{r!}\int_0^t(t-s)^r\partial_s^{r+1}u_{\adv}(s)\,\d s.
  $$
  Taking norms yields
  $$
    \norm{u_{\adv}(t)-\widetilde u_{\adv}^{(r)}(t)}
    \leq\frac{t^{r+1}}{(r+1)!}
    \sup_{0\leq s\leq t_0}\norm{\partial_s^{r+1}u_{\adv}(s)}.
  $$
  Finally,
  $$
    \begin{aligned}
      |I_{\adv}(t)-\widehat I_{\adv}(t)|\leq\bigl(\norm{u_{\adv}(t)}+\norm{\widetilde u_{\adv}^{(r)}(t)}\bigr)
      \norm{u_{\adv}(t)-\widetilde u_{\adv}^{(r)}(t)}\lesssim t^{r+1},
    \end{aligned}
  $$
  since both norms in the first factor are bounded on $[0,t_0]$ by the regularity assumption.
\end{proof}

\section{Proof of Lemma \ref{lem:tangent}}\label{sec:proof_tangent}
\begin{proof}
  Set $P=U^\top U$ and $g=g_{m,U}$. Differentiation with respect to the coefficient gives
  $$
    \left.\partial_c(cg)\right|_{c=1}=g,
  $$
  while differentiation with respect to the center yields
  $$
    \partial_{m_\alpha}g
    =
    g\,e_\alpha^\top P(x-m),
    \qquad
    \alpha=1,\dots,d.
  $$
  Since $P$ is invertible, these derivatives generate all linear modulations of $g$.

  For an upper-triangular perturbation $V$ of $U$, one has
  $$
    D_Ug[V]
    =
    -\frac12
    g\,(x-m)^\top
    \left(
    V^\top U+U^\top V
    \right)
    (x-m).
  $$
  The map $
    V
    \longmapsto
    V^\top U+U^\top V
  $
  is an isomorphism from the space of upper-triangular matrices onto $\text{Sym}(d)$. Indeed, the two spaces have the same dimension, and if
  $V^\top U+U^\top V=0$,
  then $VU^{-1}$ is both upper triangular and skew-symmetric, and therefore vanishes. Hence the shape derivatives generate all quadratic modulations of $g$. Together with the coefficient and center derivatives, they generate precisely $g\mathcal P_{\leq 2}$.
\end{proof}

\section*{Acknowledgements}

Funded by the European Union (ERC, COCOA, 101170147). Views and opinions expressed are however those of the authors only and do not necessarily reflect those of the European Union or the European Research Council. Neither the European Union nor the granting authority can be held responsible for them. The authors acknowledge the use of ChatGPT 5.6 / 6 / 6.1 Sol and ChatGPT 6 Astra for coding assistance, particularly the efficient implementation of polynomial evaluation and analytic gradient computation, as well as literature searches and manuscript proofreading.

\bibliographystyle{plain}
\bibliography{references}

@article{LasserLubich2020,
  author  = {Lasser, Caroline and Lubich, Christian},
  title   = {Computing quantum dynamics in the semiclassical regime},
  journal = {Acta Numer.},
  volume  = {29},
  pages   = {229--401},
  year    = {2020},
  doi     = {10.1017/S0962492920000033}
}

@article{AMP06,
  author  = {Arendt, Wolfgang and Metafune, Giorgio and Pallara, Diego},
  title   = {{Schr{\"o}dinger} operators with unbounded drift},
  journal = {J. Operator Theory},
  volume  = {55},
  number  = {1},
  pages   = {185--211},
  year    = {2006},
  url     = {https://www.uni-ulm.de/fileadmin/website_uni_ulm/mawi.inst.020/arendt/downloads/pubbib/short/2006-AreMetPal-SchOprUnbDrf.pdf}
}

@article{LLMR22,
  author   = {Angiuli, Luciana and Lorenzi, Luca and Mangino, Elisabetta M.
              and Rhandi, Abdelaziz},
  title    = {Generation results for vector-valued elliptic operators with
              unbounded coefficients in {$L^p$} spaces},
  journal  = {Ann. Mat. Pura Appl. (4)},
  fjournal = {Annali di Matematica Pura ed Applicata. Series IV},
  volume   = {201},
  year     = {2022},
  number   = {3},
  pages    = {1347--1379},
  doi      = {10.1007/s10231-021-01160-z}
}

@article{jahnke_lubich2000,
  author  = {Jahnke, Tobias and Lubich, Christian},
  title   = {Error Bounds for Exponential Operator Splittings},
  journal = {BIT},
  volume  = {40},
  number  = {4},
  pages   = {735--744},
  year    = {2000},
  doi     = {10.1023/A:1022396519656}
}

@article{iserles_kropielnicka2024,
  author  = {Iserles, Arieh and Kropielnicka, Karolina},
  title   = {An elementary approach to splittings of unbounded operators},
  journal = {Comput. Math. Appl.},
  volume  = {176},
  pages   = {29--34},
  year    = {2024},
  doi     = {10.1016/j.camwa.2024.08.031}
}

@incollection{more1978,
  author    = {Mor{\'e}, Jorge J.},
  title     = {The {Levenberg--Marquardt} algorithm: Implementation and theory},
  editor    = {Watson, George Alistair},
  booktitle = {Numerical Analysis},
  series    = {Lecture Notes in Mathematics},
  volume    = {630},
  pages     = {105--116},
  year      = {1978},
  publisher = {Springer},
  doi       = {10.1007/BFb0067700}
}

@article{CHANG19701,
  title   = {A practical difference scheme for {F}okker-{P}lanck equations},
  journal = {J. Comput. Phys.},
  volume  = {6},
  number  = {1},
  pages   = {1-16},
  year    = {1970},
  issn    = {0021-9991},
  doi     = {https://doi.org/10.1016/0021-9991(70)90001-X},
  author  = {Chang, Julius S. and Cooper, Gilbert}
}

@book{GMM2017,
  author    = {Gualdani, Maria Pia and Mischler, St{\'e}phane and Mouhot, Cl{\'e}ment},
  title     = {Factorization of Non-Symmetric Operators and Exponential {H}-Theorem},
  series    = {M{\'e}moires de la Soci{\'e}t{\'e} Math{\'e}matique de France},
  number    = {153},
  publisher = {Soci{\'e}t{\'e} Math{\'e}matique de France},
  year      = {2017},
  isbn      = {978-2-85629-874-9}
}

@article{wen2023couplingparameterparticledynamics,
  title   = {Coupling parameter and particle dynamics for adaptive sampling in {Neural} {Galerkin} schemes},
  volume  = {462},
  issn    = {0167-2789},
  doi     = {10.1016/j.physd.2024.134129},
  journal = {Phys. D},
  author  = {Wen, Yuxiao and Vanden-Eijnden, Eric and Peherstorfer, Benjamin},
  year    = {2024},
  pages   = {134129}
}

@article{Bruna2024,
  author   = {Bruna, Joan and Peherstorfer, Benjamin and Vanden-Eijnden, Eric},
  title    = {Neural {Galerkin} schemes with active learning for high-dimensional evolution equations},
  fjournal = {Journal of Computational Physics},
  journal  = {J. Comput. Phys.},
  issn     = {0021-9991},
  volume   = {496},
  pages    = {112588},
  year     = {2024},
  language = {English},
  doi      = {10.1016/j.jcp.2023.112588},
  zbmath   = {7807993}
}

@article{MW2018,
  title    = {On the a posteriori error analysis for linear {F}okker–{P}lanck models in convection-dominated diffusion problems},
  journal  = {Appl. Math. Comput.},
  volume   = {339},
  pages    = {779-804},
  year     = {2018},
  issn     = {0096-3003},
  doi      = {https://doi.org/10.1016/j.amc.2018.05.050},
  url      = {https://www.sciencedirect.com/science/article/pii/S0096300318304673},
  author   = {Matculevich, Svetlana and Wolfmayr, Monika}
}

@article{RHCTT2024,
  title    = {Scalable implicit solvers with dynamic mesh adaptation for a relativistic drift-kinetic {F}okker–{P}lanck–{B}oltzmann model},
  journal  = {J. Comput. Phys.},
  volume   = {507},
  pages    = {112954},
  year     = {2024},
  issn     = {0021-9991},
  doi      = {https://doi.org/10.1016/j.jcp.2024.112954},
  url      = {https://www.sciencedirect.com/science/article/pii/S0021999124002031},
  author   = {Johann Rudi and Max Heldman and Emil M. Constantinescu and Qi Tang and Xian-Zhu Tang}
}

@misc{BCM2025,
  author        = {Bon, Daan and Caris, Benjamin and Mula, Olga},
  title         = {Stable Nonlinear Dynamical Approximation with Dynamical Sampling},
  year          = {2025},
  howpublished  = {arXiv:2505.11938},
  archiveprefix = {arXiv},
  eprint        = {2505.11938},
  primaryclass  = {math.NA},
  url           = {https://arxiv.org/abs/2505.11938}
}

@article{lozinski2003fast,
  title     = {A fast solver for {F}okker--{P}lanck equation applied to viscoelastic flows calculations: {2D FENE} model},
  author    = {Lozinski, Alexei and Chauviere, C{\'e}dric},
  journal   = {J. Comput. Phys.},
  volume    = {189},
  number    = {2},
  pages     = {607--625},
  year      = {2003},
  publisher = {Elsevier}
}

@article{masud_bergman2005,
  author  = {Masud, Arif and Bergman, Lawrence A.},
  title   = {Application of Multi-Scale Finite Element Methods to the Solution of the {Fokker--Planck} Equation},
  journal = {Comput. Methods Appl. Mech. Engrg.},
  volume  = {194},
  pages   = {1513--1526},
  year    = {2005},
  doi     = {10.1016/j.cma.2004.06.041}
}

@article{liu_li_zha_zhou2022,
  author  = {Liu, Shu and Li, Wuchen and Zha, Hongyuan and Zhou, Haomin},
  title   = {Neural Parametric {Fokker--Planck} Equation},
  journal = {SIAM J. Numer. Anal.},
  volume  = {60},
  number  = {3},
  pages   = {1385--1449},
  year    = {2022},
  doi     = {10.1137/20M1344986}
}

@article{scharpenberg_lukacova2012,
  author  = {Scharpenberg, Moriz Dirk and Luk{\'a}{\v c}ov{\'a}-Medvid'ov{\'a}, Maria},
  title   = {Adaptive Gaussian Particle Method for the Solution of the {Fokker--Planck} Equation},
  journal = {ZAMM},
  volume  = {92},
  number  = {10},
  pages   = {770--781},
  year    = {2012},
  doi     = {10.1002/zamm.201100088}
}

@misc{caucheteux2026unifying,
  title        = {A Unifying View of Variational Generative Wasserstein Flows},
  author       = {Caucheteux, Paul and Bonet, Cl{\'e}ment and Korba, Anna},
  howpublished = {arXiv:2605.31369},
  year         = {2026}
}

@misc{gonnella2026dynamical,
  title        = {A Dynamical Approximation Scheme on the {S}tiefel manifold for {W}asserstein Gradient Flows},
  author       = {Gonnella, Isabella Carla and Mula, Olga and Pichi, Federico and Rozza, Gianluigi},
  howpublished = {arXiv:2607.17235},
  year         = {2026}
}

@misc{liero2025evolution,
  title        = {Evolution of {G}aussians in the {Hellinger-Kantorovich-Boltzmann} gradient flow},
  author       = {Liero, Matthias and Mielke, Alexander and Tse, Oliver and Zhu, Jia-Jie},
  howpublished = {arXiv:2504.20400},
  year         = {2025}
}

@misc{LLB2026,
  title         = {An equivalence of moment closure and nonlinear variational approximation of the Fokker-Planck equation for dilute polymeric flow},
  author        = {Caroline Lasser and Stephan B. Lunowa and Barbara Wohlmuth},
  year          = {2026},
  eprint        = {2602.04644},
  archiveprefix = {arXiv},
  primaryclass  = {math.SP},
  url           = {https://arxiv.org/abs/2602.04644}
}

@inproceedings{LambertChewiBachBonnabelRigollet2022,
  title     = {Variational Inference via {W}asserstein Gradient Flows},
  author    = {Lambert, Marc and Chewi, Sinho and
               Bach, Francis and Bonnabel, Silv{\`e}re and
               Rigollet, Philippe},
  booktitle = {Advances in Neural Information Processing Systems},
  volume    = {35},
  year      = {2022}
}

@article{chen_majda2018,
  author  = {Chen, Nan and Majda, Andrew J.},
  title   = {Efficient Statistically Accurate Algorithms for the {Fokker--Planck} Equation in Large Dimensions},
  journal = {J. Comput. Phys.},
  volume  = {354},
  pages   = {242--268},
  year    = {2018},
  doi     = {10.1016/j.jcp.2017.10.022}
}

@article{vishwajeet_singla2018,
  author  = {Vishwajeet, Kumar and Singla, Puneet},
  title   = {Adaptive Split/Merge-Based {G}aussian Mixture Model Approach for Uncertainty Propagation},
  journal = {J. Guid. Control Dyn.},
  volume  = {41},
  number  = {3},
  pages   = {603--617},
  year    = {2018},
  doi     = {10.2514/1.G002801}
}

@article{KS2009,
  author    = {Knezevic, David J. and S{\"u}li, Endre},
  title     = {Spectral {G}alerkin approximation of {Fokker-Planck} equations with unbounded drift},
  journal   = {ESAIM Math. Model. Numer. Anal.},
  year      = {2009},
  volume    = {43},
  number    = {3},
  pages     = {445--485},
  publisher = {EDP Sciences}
}

@misc{jiang_wang_zhou2026,
  author        = {Jiang, Shidong and Wang, Dong and Zhou, Qi},
  title         = {A Fast Sum-of-Gaussians Algorithm for the High-Dimensional Fractional {Fokker--Planck} Equation},
  year          = {2026},
  howpublished  = {arXiv:2606.28184},
  eprint        = {2606.28184},
  archiveprefix = {arXiv},
  primaryclass  = {math.NA},
  doi           = {10.48550/arXiv.2606.28184}
}

@incollection{Robert2004,
  author    = {Robert, Christian P.},
  title     = {Monte Carlo Methods},
  publisher = {Wiley Online Library},
  year      = {2004}
}

@inproceedings{DaumHuang2003,
  author    = {Daum, Fred and Huang, Jim},
  title     = {Curse of Dimensionality and Particle Filters},
  booktitle = {2003 IEEE Aerospace Conference Proceedings},
  volume    = {4},
  pages     = {4\_1979--4\_1993},
  year      = {2003},
  publisher = {IEEE},
  doi       = {10.1109/AERO.2003.1235126}
}

@article{CarrilloCraigPatacchini2019,
  author  = {Carrillo, Jos{\'e} A. and Craig, Katy and Patacchini, Francesco S.},
  title   = {A Blob Method for Diffusion},
  journal = {Calc. Var. Partial Differential Equations},
  volume  = {58},
  number  = {2},
  pages   = {53},
  year    = {2019},
  doi     = {10.1007/s00526-019-1486-3}
}

@article{JinMarkowichSparber2011,
  author  = {Jin, Shi and Markowich, Peter and Sparber, Christof},
  title   = {Mathematical and Computational Methods for Semiclassical
             Schr{\"o}dinger Equations},
  journal = {Acta Numer.},
  volume  = {20},
  pages   = {121--209},
  year    = {2011},
  doi     = {10.1017/S0962492911000031}
}

@article{AndersonFarazmand2024,
  author  = {Anderson, William and Farazmand, Mohammad},
  title   = {Fisher Information and Shape-Morphing Modes for Solving
             the Fokker--Planck Equation in Higher Dimensions},
  journal = {Appl. Math. Comput.},
  volume  = {467},
  pages   = {128489},
  year    = {2024},
  doi     = {10.1016/j.amc.2023.128489}
}

@article{TerejanuSinglaSinghScott2008,
  author  = {Terejanu, Gabriel and Singla, Puneet and Singh, Tarunraj
             and Scott, Peter D.},
  title   = {Uncertainty Propagation for Nonlinear Dynamic Systems
             Using Gaussian Mixture Models},
  journal = {J. Guid. Control Dyn.},
  volume  = {31},
  number  = {6},
  pages   = {1623--1633},
  year    = {2008},
  doi     = {10.2514/1.36247}
}

@inproceedings{MukherjeeRaiSinglaSinghPatra2016,
  author    = {Mukherjee, Arpan and Rai, Rahul and Singla, Puneet
               and Singh, Tarunraj and Patra, Abani},
  title     = {An Adaptive Gaussian Mixture Model Approach Based Framework
               for Solving Fokker--Planck Kolmogorov Equation Related to
               High Dimensional Dynamical Systems},
  booktitle = {Proceedings of the ASME 2016 International Design Engineering
               Technical Conferences and Computers and Information in
               Engineering Conference},
  volume    = {2B},
  pages     = {V02BT03A058},
  year      = {2016},
  publisher = {American Society of Mechanical Engineers},
  doi       = {10.1115/DETC2016-60312}
}

@article{ChenBillingsLuo1989,
  author  = {Chen, Sheng and Billings, Stephen A. and Luo, Wan},
  title   = {Orthogonal least squares methods and their application to non-linear system identification},
  journal = {International Journal of Control},
  volume  = {50},
  number  = {5},
  pages   = {1873--1896},
  year    = {1989}
}

@article{ChenCowanGrant1991,
  author  = {Chen, Sheng and Cowan, Colin F. N. and Grant, Peter M.},
  title   = {Orthogonal least squares learning algorithm for radial basis function networks},
  journal = {IEEE Transactions on Neural Networks},
  volume  = {2},
  number  = {2},
  pages   = {302--309},
  year    = {1991}
}

@article{Zhang2011,
  author  = {Zhang, Tong},
  title   = {Adaptive forward-backward greedy algorithm for learning sparse representations},
  journal = {IEEE Transactions on Information Theory},
  volume  = {57},
  number  = {7},
  pages   = {4689--4708},
  year    = {2011}
}

@inproceedings{CrouseWillettPattipatiSvensson2011,
  author    = {Crouse, David F. and Willett, Peter and Pattipati, Krishna and Svensson, Lennart},
  title     = {A look at {G}aussian mixture reduction algorithms},
  booktitle = {Proceedings of the 14th International Conference on Information Fusion},
  pages     = {1--8},
  year      = {2011}
}

@article{BeylkinMonzonYang2019,
  author  = {Beylkin, Gregory and Monz{\'o}n, Lucas and Yang, Xinshuo},
  title   = {Reduction of multivariate mixtures and its applications},
  journal = {Journal of Computational Physics},
  volume  = {383},
  pages   = {94--124},
  year    = {2019},
  doi     = {10.1016/j.jcp.2019.01.015}
}

@book{AAR,
  author    = {Andrews, George E. and Askey, Richard and Roy, Ranjan},
  title     = {Special functions},
  series    = {Encyclopedia of Mathematics and its Applications},
  volume    = {71},
  publisher = {Cambridge University Press, Cambridge},
  year      = {1999},
  pages     = {xvi+664},
  doi       = {10.1017/CBO9781107325937}
}

\end{document}